%% file: consjoa020511.tex
\documentclass {article}
\input xy
\xyoption{all}

\input xy
\xyoption{all}
\usepackage {amsmath}
\usepackage{amssymb}
\usepackage{amsthm}
\usepackage {amsfonts}
\usepackage{graphicx}
\usepackage{hyperref}
\usepackage{epsfig}

\newtheorem{thm}{Theorem}[section]
\newtheorem{prop}[thm]{Proposition}
\newtheorem{cor}[thm]{Corollary}
\newtheorem{lem}[thm]{Lemma}
\newcommand{\color}[6]{}
\usepackage{amssymb, amsmath, epsfig, graphics, amscd}

\usepackage{amsthm}
\newtheorem{rem}[thm]{Remark}

\newtheorem{defn}[thm]{Definition}

\newtheorem{defnprop}[thm]{Definition/Proposition}
\newcommand{\Cp}{\mathbb{C}}

\title {Consistency conditions for brane tilings}
\author {Ben Davison}
\date {5th May 2011}

\begin {document}
\maketitle
\begin{abstract}
Given a brane tiling on a torus, we provide a new way to prove and generalise the recent results of Szendroi, Mozgovoy and Reineke regarding the Donaldson-Thomas theory of the moduli space of framed cyclic representations of the associated algebra. Using only a natural cancellation-type consistency condition, we show that the algebras are 3-Calabi-Yau, and calculate Donaldson-Thomas type invariants of the moduli spaces. Two new ingredients to our proofs are a grading of the algebra by the path category of the associated quiver modulo relations, and a way of assigning winding numbers to pairs of paths in the lift of the brane tiling to the universal cover. These ideas allow us to generalise the above results to all consistent brane tilings on $K(\pi,1)$ surfaces. We also prove a converse: no consistent brane tiling on a sphere gives rise to a 3-Calabi-Yau algebra.
\end{abstract}
\section{Introduction}
This paper is concerned with quiver algebras along with superpotential relations arising from brane tilings.  In particular, we are interested in the question of calculating Donaldson-Thomas invariants for them, as defined in \cite{conifold}, and with the question of when such algebras can be expected to be Calabi-Yau.  The results addressing the first question mark a mild generalisation of the results of \cite{MR}.  This paper expands upon \cite{MR} in two ways: we assume weaker hypotheses, and give proofs for brane tilings on higher genus surfaces, and not just the torus.  Moreover, it is hoped that the broad idea which enables us to prove these theorems will find wider application, and lead to a better understanding of the relevant algebras.  The observation is that F-term equivalence of paths (the equivalence of paths arising from the superpotential relations) can be represented by a restriction of \textit{homotopy} equivalence of paths.  By considering chains of F-term equivalent paths, obtained by sequences of `basic equivalences,' that avoid a tile $j$ when lifted to the universal cover of the surface the brane tiling lives on, we obtain a homotopy of paths that avoids this tile, and so we can use the concept of \textit{winding} to describe this F-term equivalence class.  While this paper was being prepared, similar results were obtained independently by Nathan Broomhead in \cite{longout}.
\subsection{Acknowledgements}
I would like to thank Bal\'{a}zs Szendr\H{o}i for introducing me to most of this material, and Kazushi Ueda for his willingness to discuss ideas regarding consistency conditions.  I am indebted to Alastair King and Nathan Broomhead for sharing their thoughts, and for suggestions and corrections that led to the rewrite of the final section of this paper.
\section{Basic notions}
\begin{figure}
\centering
\input{fterm.tex}
\caption{A basic F-term equivalence takes us fom one broken loop to the other.}
\label{Fterm}
\end{figure}
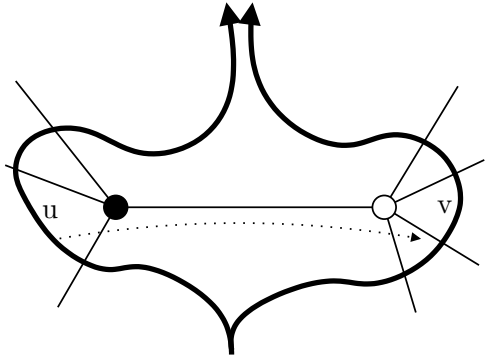

We briefly recall some definitions and fix some notation.  We will denote the compact real surface of genus $g$ by $M_g$.  A \textit{brane tiling} of genus $g$ is a tiling $F$ of $M_g$ by convex polygonal tiles and an identification of the 1-skeleton of $F$ with a bipartite graph, along with a chosen colouring of the vertices by black and white, such that each edge contains one vertex from each coloured set.  We will denote this coloured bipartite graph $\Gamma$ in general.  To such data we associate a quiver $Q$, by the following rules: the underlying graph of $Q$ is the dual graph of $F$, and the edges are directed so that they go clockwise around the black vertices.  Throughout this paper we will talk of edges of the tiling and of the quiver interchangeably - clearly they stand in a 1-1 correspondence.  To the quiver $Q$ we associate the usual path algebra $\Cp Q$.  Note that in our notation, if $u$, $v$ are paths in the path algebra, with $u$ ending at the vertex that $v$ begins at, then $uv$ is the path obtained by first going along $u$, and then $v$.  Our modules are generally right modules. \smallbreak 
We now consider the element $W\in \Cp Q/[\Cp Q, \Cp Q]$ given by taking the sum, over the vertices of $\Gamma$, of the cycle of arrows going around the vertex, signed so that those going around the black vertices occur with coefficient $+1$, and those going around the white vertices occur with coefficient $-1$.  To each edge $E$ of the quiver, there is an associated edge variable $X_E\in \Cp Q$, and there is a noncommutative derivation 
\[
\partial/\partial X_E:\Cp Q/[\Cp Q,\Cp Q]\rightarrow \Cp Q
\]
given on each monomial by cyclically permuting each instance of $X_E$ to the front and then deleting it.  We define 
\[A_F=\Cp Q/\langle\{(\partial/\partial X_E)W;E\in E(Q)\}\rangle
\]
where here and elsewhere $E(Q)$ are the edges of $Q$ and $V(Q)$ are the vertices. $A_F^+$ is the subalgebra generated by paths of length at least 1.  The relations $\{(\partial/\partial X_E)W;E\in E(Q)\}$ will be referred to as the `F-terms', borrowing notation from the physics literature.  F-term equivalence of two paths in the path algebra, then, is nothing more than equality modulo the two sided ideal generated by these relations.  The relations have form as depicted in Figure (\ref{Fterm}).\smallbreak
Let $\tilde{M_g}$ denote the universal covering space of $M_g$.  We denote the lift of $F$ (respectively $\Gamma$, $Q$) to $\tilde{M_g}$ by $\tilde{F}$ (respectively $\tilde{\Gamma}$, $\tilde{Q}$).  When we talk of paths in $\tilde{F}$ we mean paths along tiles (precisely, sequences of tiles such that consecutive elements share edges, and such that in the dual quiver the arrow goes the right way across this edge).  This is of course the same as a path in the lifted dual quiver, and such a path gives an element of $A_F$.  The tiles of $F$ or $\tilde{F}$ are labelled by the letters $i,j,k$ in this paper.  Paths are labelled $u,v,w$.  If $u$ is a path in $\tilde{F}$ or $F$, $h(u)$ will denote the tile that $u$ ends at, and $t(u)$ will denote the tile that $u$ starts at. If $i$ is a tile in $F$, the tiles in $\tilde{F}$ that are lifted from $i$ are called $i$-tiles.  Where we wish to distinguish between the tiles of $F$ and the tiles of $\tilde{F}$, we will denote the tiles of $F$ by single letters (e.g. $i$), and the tiles of $\tilde{F}$ by the single letter corresponding to the tile they project to, along with a subscript (e.g. $i_1$).  Where there is no question of ambiguity, we will use single letters for both, in order to avoid clutter.  We will write $[u]$ for the element of $A_F$ that the path $u$ represents, as well as the F-term equivalence class that $u$ belongs to.  To avoid any confusion: a letter $u$ outside of square brackets will always be a path along tiles in $\tilde{F}$ or $F$, and will not be used to denote an element of $A_F$.  In particular, we will only write $u=v$ when $u$ is equal to $v$ \textit{as a sequence of tiles}, a stronger condition than $[u]=[v]$, which merely implies that $u$ and $v$ are F-term equivalent.  We will set $u\sim v$ to mean $[u]=[v]$, i.e. $u$ and $v$ are F-term equivalent.\smallbreak
\begin{defn}
An $i$-path is a path in $\tilde{F}$ starting from some fixed $i$-tile.  An $i$-loop is an $i$-path in $\tilde{F}$ that terminates at an $i$-tile, i.e. one that projects to a loop from the tile $i$ in $M_g$.  A closed $i$-loop is an $i$-loop that terminates at the original $i$-tile in $\tilde{F}$.  A simple $i$-loop (denoted $\omega_i$) is a closed $i$-loop that consists of a single revolution around one of the vertices of the chosen $i$-tile, i.e. the completion of either of the loops in Figure (\ref{Fterm}).  For each arrow $a$ in $\tilde{F}$ there are two F-term equivalent `broken loops' from $h(a)$ to $t(a)$, which we denote $\omega^{br}_a$, when we do not wish to distinguish between them (see Figure (\ref{Fterm})).  If we do wish to distinguish between them then we write $\omega^{br}_{a,B}$ or $\omega^{br}_{a,W}$, for the broken loop going around the black (respectively white) vertex.  If $u$ is a path, $u'$ will denote the path with the last arrow deleted, $'u$ the path with the first arrow deleted.
\end{defn}
\begin{rem}
\label{transport}
It is easily shown that all simple $i$-loops are F-term equivalent, and that if $u$ is a path in $\tilde{F}$, $\omega_{t(u)}u\sim u\omega_{h(u)}$.  We will frequently use the fact that simple loops can be transported in this way along paths while maintaining F-term equivalence.
\end{rem}
\begin{defn}
A path $u$ in the tiling $\tilde{F}$ is \textit{minimal} if it is not F-term equivalent to any path that contains a simple loop.  By Remark (\ref{transport}) this is equivalent to the condition that it cannot be written $u\sim v\omega_{h(v)}$ for some $v$.
\end{defn}
We next formalize a generalisation of algebra gradings for quiver algebras.  Presumably this idea is standard - it is the natural way to define a `grading with many objects'.
\begin{defn}
Let $\mathcal{C}$ be a category, and let $A$ be an algebra over an arbitrary field $\mathbb{K}$.   A \textit{grading} of $A$ by $\mathcal{C}$ is given by a decomposition of the underlying vector space of $A$
\begin{equation*}
A\cong\bigoplus_{\phi \in mor(\mathcal{C})}V_{\phi}
\end{equation*}
where we denote the morphisms of $\mathcal{C}$ by $mor(\mathcal{C})$.  This decomposition is required to satisfy the two properties
\begin{enumerate}
\label{gradlist}
\item
If $x,y,z$ are objects of $\mathcal{C}$, and $\phi\in Hom(x,y)$ and $\psi\in Hom(y,z)$, then for all $a\in V_{\phi}$ and $b\in V_{\psi}$ $ba\in V_{\psi\phi}$.
\item
If $x,y,z,w$ are objects of $\mathcal{C}$, and $y\neq z$, and $\phi\in Hom(x,y)$ and $\psi\in Hom(z,w)$, then for all $a\in V_{\phi}$ and $b\in V_{\psi}$ $ba=0$.
\end{enumerate}
Let $A$ be a $\mathcal{C}$-graded algebra.  A $\mathcal{C}$-graded right module $M$ over $A$ is given by a decomposition
\begin{equation*}
M\cong\bigoplus_{\phi \in mor(\mathcal{C})}W_{\phi}
\end{equation*}
such that the composition $M\otimes A\rightarrow M$ satisfies conditions analogous to (1) and (2) above.
\end{defn}
This can be generalised to a grading of a category by another category.  The canonical example is the grading of $\mathbb{C}Q$ by the unlinearised path category of $Q$, for an arbitrary quiver $Q$.  
\smallbreak
We next come to a definition of consistency.  This is equivalent to Condition 4.12 of \cite{MR}.
\begin{defn}
\label{consistency}
Two paths $u$ and $v$ ending at the tile $j$ are \textit{weakly equivalent} if there is some $n\in\mathbb{Z}_{\geq 0}$ such that $u\omega_j^n\sim v\omega_j^n$.  A tiling is \textit{consistent} if weak equivalence implies F-term equivalence.
\end{defn}
This is essentially a cancellation law.  It is equivalent to the condition that for all paths $u$, $v$, $w$, if $[u][w]\neq 0$ and $[u][w]=[v][w]$ then $[u]=[v]$.  The set of paths in the quiver (with relations) naturally has the structure of a semigroupoid $G_0$, which is in fact a category, if we include zero length paths at the tiles $i$ of $F$, i.e. $u_i$, such that $[u_i]=e_i$, where $e_i$ is the idempotent corresponding to the tile $i$.  Precisely:
\begin{defn}
$G_0$ is the category having, as objects, the tiles of $F$.  Morphisms in $G_0$ are given by sequences of tiles $(i_0,i_1,\ldots,i_n)$ in $F$, where $n\geq 0$, up to F-term equivalence, that correspond to paths in the dual quiver $Q$ that respect the orientation of $Q$.  Composition of morphisms is given by composition of paths.  The category $G_1$ is given by formally inverting morphisms of $G_0$.
\end{defn}
In \cite{MR} it is shown that under certain circumstances $A_F$ possesses a natural torus action that induces a grading.  We replace this with the natural grading of $A_F$ by the category $G_1$, induced by the natural functor from $G_0$ to $G_1$.  The following is a useful technical lemma:
\begin{lem}
\label{loopylemma}
If $u$ is a closed $i$-loop in an arbitrary brane tiling $F$, then there exist $n,m\in \mathbb{Z}_{\geq 0}$ such that $u \omega_{i}^n \sim \omega_{i}^m$.
\end{lem}
\begin{proof}  Let $v$ be a path along the tiles of $\tilde{F}$, where, for now, we don't assume that $v$ passes through the edges in the correct direction.  To such a path we associate a $u$ that does go through every edge in the right direction by the following rule: we go along the edges of $v$, and each time an edge $e$ passes in the correct direction between the tiles of $\tilde{F}$ we add that edge to $u$, and each time an edge passes the opposite way, we add one of the broken loops associated to that edge to $u$.  Given such a path $v$ we claim that there is a $n\in \mathbb{N}$ such that, if $u$ is the path associated to $v$ by following the given rule, $u \omega_{i}^n\sim \omega_{i}^m$ for some $m\in\mathbb{N}$.  This clearly suffices to prove the lemma.  Let $v$ be such a path.  We can assume that $v$ visits no tile twice, for otherwise we may prove the claim after breaking $v$ into two disjoint loops.  Note that adding an edge, and then its inverse, results in adding a simple loop under the rule for associating legal paths to paths that pass the wrong way through some edges.  The claim is easy to prove in the case where the number of vertices in the interior of the loop $v$ is 1.  But we can easily reduce to this case by adding extra paths along the interior of the loop along with their inverses, and we're done.
\end{proof}

In fact the proof clearly gives us more.  As in \cite{MR} we deduce that
\begin{lem}
\label{toruscone}
If $F$ is a consistent tiling, and $u$ and $v$ have the same startpoints and endpoints in $\tilde{F}$, they satisfy $u\sim v\omega_{t(v)}^n$ for some $n\in\mathbb{Z}_{\geq 0}$, or $v\sim u\omega_{t(u)}^n$, for some $n\in\mathbb{Z}_{\geq 0}$.
\end{lem}
\begin{lem}
\label{minimalexist}
Let $F$ be a consistent tiling.  Let $i,j\in \tilde{F}$ be two tiles.  Then there is a minimal path $u$ that starts at $i$ and ends at $j$.
\end{lem}
\begin{proof}
We prove this by induction on $n\geq 0$, the shortest length of a path from $i$ to $j$.  If there is a path of length $0$ from $i$ to $j$ we deduce that $i=j$ and so the trivial path of length zero, starting at $i$, is a minimal path from $i$ to $j$.  So assume $i\neq j$, and, using the inductive hypothesis, assume there is a $k\in \tilde{F}$ such that there is an arrow $a$ from $k$ to $j$ in $\tilde{Q}$, the lift of the dual quiver for $F$, and a minimal path $u$ starting at $i$ and ending at $k$.  Assume also that $ua$ is not minimal, for else we deduce immediately that there is a minimal path from $i$ to $j$.  It follows that $ua\sim v\omega_j$ for some path $v$.  We can pick $\omega_j$ to go through $k$, since all simple loops are F-term equivalent.  This gives us a decomposition $ua\sim v\omega_{a}^{br} a$.  By consistency we deduce that $u\sim v\omega_{a}^{br}$.  Now from the definition of minimality of paths we deduce that a sub-path of a minimal path is minimal, from which we obtain that $v$ is a minimal path to $j$.
\end{proof}
\begin{rem}
\label{funnyloops}
If $F$ is a consistent tiling then, using Lemma (\ref{toruscone}) and Remark (\ref{transport}), we deduce that a path $u$ in $\tilde{F}$ is minimal if and only if it is not F-term equivalent to a path that visits some tile twice.
\end{rem}
The consistency condition can be restated as saying that the natural functor from $G_0$ to $G_1$ is faithful.  There is a natural grading of the algebra $A_F$ by the groupoid $G_1$.  Finally, it is seen that weak equivalence of two paths $u$ and $v$ is equivalent to $[u]$ and $[v]$ having the same $G_1$-grade, and consistency is equivalent to each $G_1$-graded piece of $A_F$ being 1-dimensional or 0-dimensional.  
We will make use of the following:
\begin{defnprop}
\label{ordering}
Let $F$ be a consistent brane tiling.  For $i$ a fixed tile in $\tilde{F}$, we define the partial ordering on tiles in $\tilde{F}$, denoted $\leq_i$, by the rule: $j\leq_i k$ if there is a minimal path from $i$ to $k$ passing through $j$.  We write $j<_i k$ to mean that $j\leq_i k$ and $j\neq k$.
\end{defnprop}
We prove that $\leq_i$ is a partial order:
\begin{enumerate}
\item{Transitivity:}
Let $i,j,k,l$ be tiles of $\tilde{F}$, satisfying $j \leq_i k$ and $k \leq_i l$.  Then there is a minimal path $u$ from $i$ to $l$ passing through $k$.  We decompose $u$ as $u=vw$ where $v$ is a path from $i$ to $k$, and $w$ is a path from $k$ to $l$.  Now $v$ is minimal, for otherwise $u$ certainly is not.  From Lemma (\ref{toruscone}) we deduce that all minimal paths from $i$ to $k$ are F-term equivalent, so since $j \leq_i k$ we see that $v$ is F-term equivalent to a path passing through $j$.  It follows that $u$ is F-term equivalent to a path passing through $j$.
\item{Antisymmetry:}
Let $i,j,k$ be tiles of $\tilde{F}$.  Assume that $j\leq_i k$ and $k\leq_i j$.  Then there is a minimal path $u$ from $i$ to $k$ passing through $j$.  We decompose this path as $u=vw$ where $v$ is a path from $i$ to $j$ and $w$ is a path from $j$ to $k$.  As in the check of transitivity, we deduce that $v$ is minimal.  Since we have $k \leq_i j$, we deduce, as before, that $v$ is F-term equivalent to a path passing through $k$.  It follows that if $w$ has length greater than zero then $u$ is F-term equivalent to a path that visits the tile $k$ twice, and is not minimal, by Remark (\ref{funnyloops}).  So $w$ has length zero and $j=k$.
\item{Reflexivity:}
Let $i$, $j$ be tiles of $F$.  Let $u$ be a minimal path from $i$ to $j$, which exists by Lemma (\ref{minimalexist}).  Then this path demonstrates that $j\leq_i j$. 
\end{enumerate}

\section{F-term equivalence and homotopy equivalence}
\label{homotopies}
We introduce here a vital heuristic for the proofs in this paper.  Consider again the terms coming from the noncommutative derivations of the superpotential.  They give equivalences between the `broken loops' as depicted in Figure (\ref{Fterm}).  In fact the relations are given precisely by identifying each set of broken loops.  We say there is a \textit{basic} F-term equivalence between two paths $u$ and $v$ of $\tilde{F}$ if
\[
\{u,v\}=\{w_1\omega^{br}_{a,B} w_2,w_1\omega^{br}_{a,W}w_2\}\] for some paths $w_1$, $w_2$ in $\tilde{F}$, and some arrow $a$ in $\tilde{F}$. 
\smallbreak
Let $\gamma$ be a (suitably well behaved) continuous path through $\tilde{M}_g$, in the sense of a continuous function 
\[
\gamma:[0,1]\rightarrow \tilde{M}_g.
\]
As long as $\gamma$ passes edges of $\tilde{F}$ the `right' way (i.e. with the black vertex to the right), we can naturally associate to it a path (in the sense used throughout the rest of this paper) along the dual quiver $\tilde{Q}$.  In this way we obtain a function
\begin{equation}
\label{psidef}
\Psi:Cont_F([0,1],\tilde{M}_g)\rightarrow paths_{\tilde{F}}
\end{equation}
where $Cont_F([0,1],\tilde{M}_g)$ is just the set of continuous paths such that the above construction gives a well-defined path in $\tilde{F}$, crossing each edge of $\tilde{F}$ with the right orientation, and $paths_{\tilde{F}}$ is the set of paths of $\tilde{F}$.  Note that this function is surjective.
\smallbreak
Let $u$ and $v$ be two paths of $\tilde{F}$ with the same starting point and the same endpoint.  If $u\sim v$ then there is a chain of basic F-term equivalences taking $u$ to $w$.  Let $\gamma_u$ and $\gamma_w$ be two continuous paths in $Cont_F([0,1],\tilde{M}_g)$ such that
\begin{eqnarray*}
\Psi(\gamma_u)&=&u\\
\Psi(\gamma_w)&=&w.
\end{eqnarray*}
Say $v$ is obtained from $u$ by applying a basic F-term equivalence.  Then there exist paths $u_1$ and $u_2$ in $\tilde{F}$ such that $u=u_1\omega^{br}_{a,B} u_2$ and $v=u_1\omega^{br}_{a,W} u_2$.  There exist $a,b \in [0,1]$ such that the copy of $\omega^{br}_{a,B}$ indicated in the given decomposition of $u$ is obtained by
\[
\omega^{br}_{a,B}=\Psi(\gamma_u|_{[a,b]}\circ s)
\]
where $s$ is just a scaling function $s(x)=x(b-a)+a$.  As indicated in Figure (\ref{Fterm}) there is a homotopy
\[
H:[0,1]\times[0,1]\rightarrow\tilde{M}_g
\]
satisfying $H|_{\{0\}\times[0,1]}=\gamma_u|_{[a,b]}\circ s$ and $H|_{\{1\}\times [0,1]}=\nu$, where 
\[
\nu:[0,1]\rightarrow \tilde{M}_g.
\]
is a new path satisfying $\Psi(\nu)=\omega^{br}_{a,W}$.  We pick this homotopy to be invariant along the line $[0,1]\times\{0\}$, and invariant along the line $[0,1]\times \{1\}$.  This is just the same as stipulating that the homotopy fixes the endpoints, so we can construct a new homotopy $H'$ given by
\begin{eqnarray*}
H'(t,x)&=\gamma_u(x) &\hbox{if } t\notin[a,b]\\
&=H(t,s^{-1}(x))&\hbox{otherwise.}
\end{eqnarray*}
This gives a homotopy such that $\Psi(H'|_{\{0\}\times [0,1]})=u$ and $\Psi(H'|_{\{1\}\times[0,1]})=v$.  Note that, as indicated in the figure, we can pick $H'$ so that $image(H':[0,1]^2\rightarrow \tilde{M}_g)$ does not intersect any of the tiles not included in $u$ or $v$.
\smallbreak
The crucial observation is that if there is a chain of basic F-term equivalent paths between two paths $u$ and $w$ through $\tilde{F}$, none of which go through some tile $k$ of $\tilde{F}$, then we can produce $\gamma_w,\gamma_u\in Cont_F([0,1],\tilde{M}_g)$ such that $\Psi(\gamma_w)=w$, $\Psi(\gamma_v)=v$, and an endpoint-preserving homotopy $H$ from $\gamma_w$ to $\gamma_v$ such that
\[
Image(H:[0,1]^2\rightarrow \tilde{M}_g)\cap k=\emptyset.
\]
This motivates the following definition.
\begin{defn}
A pair of F-term equivalent paths $(u,v)$ of $\tilde{F}$ is \textit{$k$-unwound} if there is a sequence of paths $u=u_0,\ldots,u_n=v$ such that $u_s$ is related to $u_{s+1}$ by a basic F-term equivalence, for every $s$, and none of the $u_s$ pass through tile $k$.  
\end{defn}
In particular, if $u,v\in paths_{\tilde{F}}$ are a pair of $k$-unwound paths then we can find a pair of paths $\gamma_u,\gamma_v\in Cont([0,1],\tilde{M}_g)$ that avoid $k$, and such that $\Psi(\gamma_u)=u$ and $\Psi(\gamma_v)=v$, and such that there is an endpoint preserving homotopy between $\gamma_u$ and $\gamma_v$, considered as paths in $\tilde{M}_g\backslash k$.  While this definition makes sense for tilings on all universal covers of oriented surfaces $\tilde{M}_g$, for $g=0$ it is not useful, since all paths between two chosen tiles in $\tilde{M}_g\backslash k$ are homotopy equivalent.  We could of course define the $K$-winding number of the ordered pair $(u,v)$, or the $H_1(\tilde{M_g}\backslash K)$-class of the ordered pair $(u,v)$ for $K$ any union of tiles of $\tilde{F}$, and $u$, $v$ two paths avoiding $K$.  However, we have no need for these added complications here (except perhaps for Remark (\ref{sphere})).\medbreak
\section{3-Calabi-Yau property}
\begin{figure}[htbp]
\begin{center}
\caption{Winding around tile $i_0$.}
\label{windup} 
\input{windup2.tex}
 
\end{center}
\end{figure}
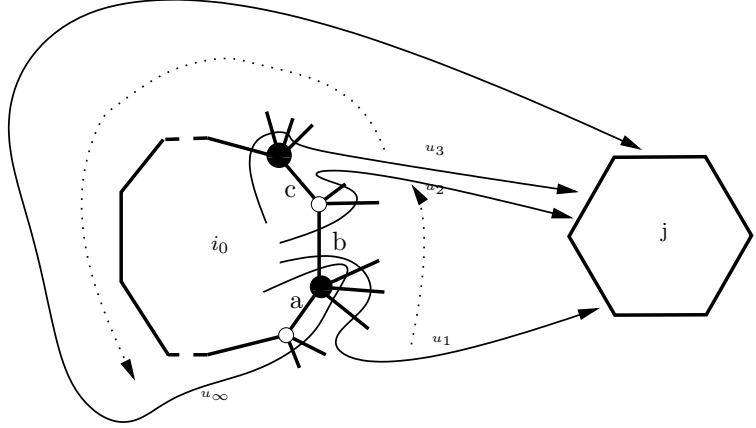

We proceed to give a detailed proof of the following theorem, in which a winding argument plays a key part.\medbreak
\begin{thm}
\label{resolution}
Let $A:=A_F$ be the algebra associated to a consistent tiling~$F$ on a surface $M_g$, for $g \geq 1$.  For each $i$ a tile of $F$ there is a $G_1$-graded projective resolution of 
\[S_i:=e_iA/e_iA^+
\]
given by 
\begin{equation}
\label{res1}
\xymatrix{
0\ar[r]&[\omega_i] A\ar[r]^(0.35){\tau}&\bigoplus_{a\in E_1}[\omega^{br}_a] A\ar[r]^{\sigma}&\bigoplus_{b\in E_2} [b] A\ar[r]^(.6){\rho}&e_i A
}
\end{equation}
\noindent
where $E_1$ is the set of edges going into $i$, and $E_2$ is the set of edges going out of $i$.  The maps are defined as follows: $\rho$ is the sum of the obvious inclusion maps.  Let $a\in E_1$.  Then there are natural inclusions $\iota_1:[\omega^{br}_{a}]A\rightarrow [b_1] A$ and $\iota_2:[\omega^{br}_{a}]A\rightarrow [b_2] A$, where $b_1$ is the first arrow in $\omega^{br}_{a,B}$ and $b_2$ is the first arrow in $\omega^{br}_{a, W}$.  We let $\sigma$ be the direct sum of $\iota_1-\iota_2$ for $a\in E_1$.  Finally, there are natural inclusions $[\omega_i] A\rightarrow[\omega^{br}_a]A$, for each $a\in E_1$, and we let $\tau$ be the sum of these.
\end{thm}
\begin{proof}
That (\ref{res1}) is a chain complex follows straight from the definitions of the maps.  All the maps in (\ref{res1}) are clearly $G_0$-graded and therefore $G_1$-graded.  There is a morphism $f$ of complexes from (\ref{res1}) to the complex consisting of zeroes in every position except for a copy of $S_i$ in the zeroeth position, induced by the map
\[
e_i A\rightarrow S_i.
\]
We must check that this induces isomorphisms on homology, i.e. that $H^i(f)$ is an isomorphism for all $i$.\smallbreak
$H^0(f)$ is an isomorphism, since the kernel of the map from $e_i A$ to $S_i$ is generated by paths of length $1$, all of which lie in the image of $\rho$.  Now we need just to check that the complex (\ref{res1}) is exact away from the zeroeth position.  Injectivity of $\tau$ follows from the injectivity of each of its components.   Now restrict attention to the $G_1$-graded piece of this chain complex corresponding to $[u]\in A$ where $u$ is a path from a chosen $i$-tile $i_0$ to some $j$-tile.  This path is fixed for the rest of the proof.  Note that for each of the summands $[v]A$, in each of the modules of (\ref{res1}), consistency implies that $\mathrm{dim}(([v]A)_{([u])})=1$ if $u$ is F-term equivalent to some path starting with $v$, and $\mathrm{dim}(([v]A)_{([u])})=0$ otherwise.\smallbreak

Define $E_2'\subset E_2$ to be the subset of $E_2$ consisting of those arrows $a$ such that there is a path $v$ with $v\sim u$ such that $v$ passes through $h(a)$.  Note that the arrow $a$ is the unique minimal path F-term equivalent to $a$, since we cannot apply any F-term equivalences to it.  Say $a\in E_2'$.  Then there is some path $v$, F-term equivalent to $u$, such that $v$ passes through $h(a)$.  By Lemma \ref{toruscone}, we can choose $v$ so that it decomposes as $v=v_1\omega_{h(a)}^n v_2$ where $v_1$ is a minimal path from $t(a)$ to $h(a)$, $v_2$ is a path from $h(a)$ to $h(u)$, and $n\geq 0$.  We deduce that $a\in E_2'$ if and only if there is a path $v$, F-term equivalent to $u$, such that $v$ starts with $a$. \medbreak
We define $E_1'$ to be the subset of $E_1$ consisting of those $a$ such that there is a path $v$ such that $v\sim u$ and $v$ passes through $t(a)$.  Again, by the fact that broken loops are the unique minimal paths F-term equivalent to broken loops, we deduce that $a\in E_1'$ if and only if there is some $v$ such that $v\sim u$ and $v$ begins with $\omega^{br}_{a}$.  \medbreak
Say, to start with, that $u$ is not a minimal path.  It follows that $u$ is F-term equivalent to some path starting with a loop $\omega_i$.  It follows that $E_1'=E_1$ and $E_2'=E_2$.  We deduce that the $[u]$-graded piece of the complex is given by
\begin{equation}
\label{res2}
\xymatrix{
0\ar[r]&\Cp\ar[r]^(0.3){\tau_{[u]}}&\bigoplus_{a\in E_1}\Cp_a\ar[r]^{\sigma_{[u]}}&\bigoplus_{b\in E_2}\Cp_b\ar[r]^(0.6){\rho_{[u]}}&\Cp
}
\end{equation}
where $\Cp_a\cong\Cp$ for each $a\in E_1\bigcup E_2$, $\tau_{[u]}$ is the diagonal map, $\rho_{[u]}$ takes vectors to the sum of their entries, and the map $\sigma_{[u]}$ takes $1_a$ to $1_{b_1}-1_{b_2}$, where $b_1$ is the edge bordering $a$ at a black vertex, and $b_2$ is the edge bordering $a$ at a white vertex.  Exactness of (\ref{res2}) is clear.
\smallbreak

So assume that $u$ is minimal.   We start by investigating what $E_1'\bigcup E_2'$ looks like (considered as edges in $\tilde{\Gamma}$, rather than arrows in the dual quiver, for now).  We say a path $w$ demonstrates that $a\in E_1'$ if $w\sim u$ and $w$ starts with $\omega^{br}_{a}$.  Similarly, a path $w$ demonstrates that $a\in E_2'$ if it starts with $a$, and $w\sim~u$.  Trivially, if a path demonstrates that an edge is in $E_1'$, it demonstrates that a neighbouring edge is in $E_2'$.  Given two edges $e$ and $f$ in $E_1'\bigcup E_2'$, we consider two paths $w_e$ and $w_f$ demonstrating that $e$ and $f$, respectively, belong to $E_1' \bigcup E_2'$.  Since $w_e$ and $w_f$ are F-term equivalent there is a sequence of paths $w_e=w_1,\ldots,w_t=w_f$ such that each $w_s$ is related to its predecessor by a basic F-term equivalence.  Each of these paths demonstrates that one or two connected members of $E_1\bigcup E_2$ is a member of $E_1'\bigcup E_2'$ (two in the case that $w_s$ starts with a broken loop).  The set of members of $E_1\bigcup E_2$ that $w_s$ demonstrates is in $E_1'\bigcup E_2'$ is connected to the set that $w_{s-1}$ demonstrates is in $E_1'\bigcup E_2'$.  We deduce that $E_1'\bigcup E_2'$ is connected.  Also, since $\omega^{br}_{a,B}\sim\omega^{br}_{a,W}$, we deduce that if $e\in E_1'$, the two edges in $E_2$ connected to $e$ are in $E_2'$.  We deduce that either

\begin{enumerate}
\item
the set of edges in $E_1'\bigcup E_2'$ forms a connected contractible set, with edges from $E_2'$ at either end, or
\item
\label{badtime}
these edges form a loop, i.e. $E_1'\bigcup E_2'=E_1\bigcup E_2$.
\end{enumerate}
Restricting to the $[u]$-grade of the complex, it becomes
\begin{equation}
\label{res3}
\xymatrix{
0\ar[r]&\bigoplus_{a\in E_1'}\Cp_a\ar[r]^{\sigma_{[u]}}&\bigoplus_{b\in E_2'}\Cp_b\ar[r]^(0.6){\rho_{[u]}}&\Cp
}
\end{equation}
where $\rho_{[u]}$ is still the summation map.  The zero at the first place is the result of our picking $u$ to be minimal, which is equivalent to $([\omega_i]A)_{[u]}=0$.  From this chain complex and our description of $E_1'\bigcup E_2'$ it is clear that exactness of (\ref{res3}) is equivalent to situation (B) not occuring.\medbreak
So say we are in situation (B) for our grade $[u]$.  We construct a sequence of basic F-term equivalences winding around tile $i_0$.  Let $a,b,c$ be edges of tile $i_0$ as indiciated in Figure (\ref{windup}).  Since we are in sitaution (B), there is a path $u_{1}$, as indicated in Figure (\ref{windup}) which is F-term equivalent to $u$, and starts with $\omega^{br}_{a,B}$.  Again, by the supposition that $E'_1\bigcup E'_2=E_1\bigcup E_2$, there is a path $u_{2}$, F-term equivalent to $u$, starting with $\omega^{br}_{c,W}$.  Both $u_1$ and $u_2$ start with the arrow $b$, and by consistency we have $'u_1\sim 'u_2$ (recall that for a path $v$ we denote by $'v$ the path obtained by deleting the first arrow of $v$).  Furthermore, $'u_1$ and $'u_2$ are $i_0$-unwound, by minimality.\medbreak
There is a basic F-term equivalence that replaces $\omega^{br}_{c,W}$ with $\omega^{br}_{c,B}$, taking $u_2$ to a new path $u_3$.  We next repeat the reasoning of the above paragraph, winding anticlockwise around the tile $i_0$.

\medbreak
We proceed in this way around tile $i_0$ until we get to a new path $u_{\infty}$ starting with $\omega^{br}_{a,B}$, F-term equivalent to $u$.  We note that none of the paths in this chain of F-term equivalent paths passes through tile $i_0$, except to start at it, by minimality of $u$.  Again, by minimality of $u$, and cancellation for consistent tilings, $'u_{\infty}$ and $'u_{1}$ are F-term equivalent and $i_0$-unwound, but this is absurd, since we have constructed $u_{\infty}$ to wind around $i_0$.
\end{proof}

\smallbreak
The following definition is taken from \cite{ginz}.
\begin{defn}
An algebra A is Calabi-Yau of dimension $d$ provided there is an A-bimodule quasi-isomorphism
\[
\xymatrix{
f:A\ar[r]^{~}&A^![d]}
\]
such that $f=f^{!}[d]$, where for an arbitrary A-bimodule the functor $M\mapsto M^!$ is given by $\mathrm{\textbf{R}Hom}_{A\otimes A^{op}\mathrm{-}mod}(M,A\otimes_{\Cp} A)$.
\end{defn}
\begin{thm}
Let $A_F$ be the algebra associated to a consistent brane tiling $F$ with dual quiver $Q$.  The complex in (\ref{res1}) is a projective resolution of $S_i$, for each $i$, if and only if $A_F$ is a Calabi-Yau algebra of dimension 3.
\label{exactCY}
\end{thm}
\begin{proof}
The method of proof is a modified version of the argument of \cite{MR}.\smallbreak
Define the algebra
\[
R:=\bigoplus_{j\in V(Q)}\mathbb{C},
\]
the direct sum of $|V(Q)|$ copies of the algebra $\mathbb{C}$.  As a vector space this is generated by elements $1_j$ for $j\in V(Q)$.  Then $A_F$ is naturally an $R$-bimodule - there is an algebra morphism $\iota:R\rightarrow A_F$ sending $1_j$ to $e_j$.  If $S$ is an $R$-bimodule, $S$ defines an endofunctor for the category of left $R$-modules via the tensor product over $R$.  We say that $S$ is (left) flat if this functor takes exact sequences to exact sequences, and is (left) faithfully flat if the exactness of $S\otimes M^{\bullet}$ for some complex $M^{\bullet}$ of left $R$-modules implies the exactness of $M^{\bullet}$.  The following facts are easy to verify.
\begin{enumerate}
\item
All $R$-bimodules are (left) flat.
\item
An $R$-bimodule $S$ is (left) faithfully flat if and only if $S\cdot 1_j\neq 0$ for all $j$.
\end{enumerate}
To each edge $a\in E(Q)$ we associate an $R$-bimodule $E_a$ in the natural way: $E_a$ has underlying vector space $\mathbb{C}$, and
\begin{eqnarray*}
1_k\cdot m\cdot1_l=\begin{cases}m&\text{if }k=t(a)\text{ and } l=h(a),\\0&\text{otherwise.}\end{cases}
\end{eqnarray*}

We pick a generator of the underlying vector space of $E_a$ and denote it $1_a$.  We define
\[
E:=\bigoplus_{a\in E(Q)} E_a.
\]
There is a natural basis for the underlying vector space of $E$ given by the elements $1_a$, for $a\in E(Q)$.\smallbreak
$R$ is a sub-algebra of $A_F$ and inherits a natural $G_1$-grading.  Furthermore we give $E$ a $G_1$-grading by setting the grade of $1_a$ to be the grade of $[a]$, considered as an element of $A_F$.
\smallbreak
Next, to each edge $a\in E(Q)$ we associate an $R$-bimodule $E^*_a$.  This has underlying vector space $\mathbb{C}$, and 
 \begin{eqnarray*}
1_k\cdot m\cdot 1_l&=\begin{cases}m&\text{if }k=h(a)\text{ and } l=t(a),\\
0&\text{otherwise.}\end{cases}
\end{eqnarray*}
We define
\[
E^*:=\bigoplus_{a\in E(Q)} E_a^*.
\]

There is a natural basis for the underlying vector space of $E^*$ given by elements denoted $1_a^*$, where $1_a^*\in E_a^*$.  This $R$-bimodule has a $G_1$-grading given by setting the grade of $1_a^*$ to be given by the path $\omega_a^{br}$.\smallbreak
Finally we define a $G_1$ graded $R$-bimodule $L$.  As an \textit{ungraded} $R$-bimodule this is isomorphic to $R$.  We shift the grading by setting the grade of the copy of $1_i$ in $L$, for some vertex $i\in V(Q)$, to be given by the path $\omega_i$.\medbreak
We must now show that the sequence
\begin{equation}
\label{ginzcomplex}
\xymatrix{
0\ar[r]&A\otimes_R L\otimes_R A\ar[r]^j& A\otimes_R E^*\otimes_R A\ar[r]& A\otimes_R E\otimes_R  A\ar[r]& A\otimes_R A\ar[r]&A
}
\end{equation}
of Proposition 5.1.9 of \cite{ginz} is exact  (where we define $A:=A_F$), if and only if the complexes of (\ref{res1}) are quasi-isomorphic to $S_i$, for each $i$.  It is known that the sequence (\ref{ginzcomplex}) is exact at every place, except possibly at the first two nonzero places.  The injectivity of the map $j$ is proved in \cite{MR}, using consistency.  From the gradings we have endowed on the bimodules $E$, $E^*$ and $L$, this is a complex of $G_1$-graded $A$-bimodules.  Denote the complex (\ref{ginzcomplex}) by $L^{\bullet}$.  Since there is only one position in which $L^{\bullet}$ may have nontrivial homology, we deduce that
\[
\sum_{n\in\mathbb{Z}} (-1)^n dim(H^{n}(L^{\bullet})_{\nu})=0
\]
for all morphisms $\nu$ in the groupoid $G_1$, if and only if $L^{\bullet}$ is exact.  This is equivalent to
\[
\sum_{n\in\mathbb{Z}} (-1)^n dim((L^{n})_{\nu})=0
\]
for all morphisms $\nu$ in the groupoid $G_1$.  For each vertex $i$ we have a complex $M^{\bullet}_i$ of right $A$-modules given by the natural morphism from the complex of (\ref{res1}) to the simple module $S_i$.  $M^{\bullet}_i$ is exact if and only if the complex (\ref{res1}) is a projective resolution of $S_i$.  We have seen, in the proof of Theorem (\ref{resolution}), that $(M^{\bullet}_i)_{\nu}$ is exact if $\nu$ is not minimal.  Furthermore, in the case in which $\nu$ is minimal, we showed that 
\[
\sum_{n\in\mathbb{Z}}(-1)^n dim((M^{n})_{\nu})=1\neq 0
\]
if $(M^{n})_{\nu}$ is not exact: this is situation (B) of that proof - which can occur only if the genus of $M_g$ is zero.  In particular, $M^{\bullet}$ can only have nontrivial homology in one position.
\smallbreak
Given a right $A_F$-module $M$ we can consider $M$ as a right $R$-module via the inclusion map $\iota:R\rightarrow A_F$.  If $M$ is a right $A_F$-module of the form $[\nu]A$, we give $M$ a \textit{left} $R$-action by setting
\begin{eqnarray*}
1_i\cdot m=\begin{cases}m&\text{if }i=t(\nu)\\
0&\text{ otherwise.}\end{cases}
\end{eqnarray*}
Then each $M^{\bullet}_i$ is a complex of $R$-bimodules.  Set
\[
M^{\bullet}:=\bigoplus_{i\in V(Q)} M^{\bullet}_i.
\]
This, again, is a $G_1$-graded complex of $R$-bimodules.  Now observe that the $R$-bimodules in the complex $A_F\otimes_R M^{\bullet}$ are the same as those of the complex $L^{\bullet}$ (though the morphisms are different).  By our criterion for (left) faithful flatness of an $R$-bimodule, we deduce that $A_F\otimes_R M^{\bullet}$ is exact if and only if all of the $M^{\bullet}_i$ are, and also $A_F\otimes_R M^{\bullet}$ can only have nontrivial homology in one position, since $M^{\bullet}$ can only have nontrivial homology in one position.  It follows that
\[
\sum_{n\in\mathbb{Z}}(-1)^n dim(H^{n}(L^{\bullet})_{\nu})=\sum_{n\in\mathbb{Z}}(-1)^n dim((L^n)_{\nu})
\]
is zero if and only if $(M^{\bullet}_{t(u)})_{\nu}$ is exact.  The theorem follows.
\end{proof}

\begin{cor}
\label{CY3}
If $A_F$ is the algebra associated to a consistent brane tiling $F$ on a surface $M_g$, with $g\geq 1$, then $A_F$ is Calabi-Yau of dimension 3.
\end{cor}

The next theorem is stated without proof.  The proof uses only those ideas found in the proof of Theorem (\ref{resolution}).  A winding argument is again used.
\begin{thm}
\label{injresolution}
If $A=A_F$ is the algebra associated to a consistent tiling $F$ on a surface of genus at least one, there is an injective resolution of $A$ as a left $A$-module:
\[
\xymatrix{
A\ar@{^{(}->}^(.25){ev}[r]&\bigoplus_{i\in V(Q)} \mathrm{Hom}_{\Cp}(e_i A, e_i A) \ar[r]^(.45){d}&\bigoplus_{a\in E(Q)} \mathrm{Hom}_{\Cp}(e_{h(a)}A,e_{t(a)}A)\ar[r]^(.85)r&\\ \ar[r]^(.2)r&\bigoplus_{a\in E(Q)} \mathrm{Hom}_{\Cp}(e_{t(a)}A,e_{h(a)}A)\ar@{>>}[r]^(.5)p&\bigoplus_{i\in V(Q)}\mathrm{Hom}_{\Cp}(e_iA,e_iA)
}
\]
and an analogous injective resolution of right $A$-modules.  Here $Q$ is the dual quiver to $F$, $V(Q)$ are its vertices, and $E(Q)$ are its edges.
\end{thm}
The maps are as follows: $ev(x)$ is the linear map given by right multiplication by $x$. For an arbitrary $a\in E(Q)$ and an element $\phi\in\bigoplus_{i\in V(Q)}\mathrm{Hom}_{\Cp}(e_iA,e_iA)$ there are two linear maps $e_{h(a)}A\rightarrow e_{t(a)}A$ given by the two routes in the diagram below:
\[
\xymatrix{
e_{h(a)}A\ar[d]^{\phi_{h(a)}}\ar[r]^{a\cdot}&e_{t(a)}A\ar[d]^{\phi_{t(a)}}\\
e_{h(a)}A\ar[r]^{a\cdot}&e_{t(a)}A
}
\]
we set
\[
(d(\phi))_a(v)=a\cdot \phi_{h(a)}(v)-\phi_{t(a)}(a\cdot v).
\]
Now let $\psi\in \bigoplus_{a\in E(Q)}\mathrm{Hom}_{\Cp}(e_{h(a)}A,e_{t(a)}A)$.  For a given $a\in E(Q)$, $\omega^{br}_{a,B}$ and $\omega^{br}_{a,W}$ give a composition of the linear maps defined by $\psi$, which we denote $\alpha_{\psi,a,B}$ and $\alpha_{\psi,a,W}$, in $\bigoplus_{a\in E(Q)}\mathrm{Hom}_{\Cp}(e_{t(a)}A,e_{h(a)}A)$.  Let $r(\psi)_{a}=\alpha_{\psi,a,B}-\alpha_{\psi,a,W}$.\smallbreak Finally, given a $\nu\in\bigoplus_{a\in E(Q)}\mathrm{Hom}_{\Cp}(e_{t(a)}A,e_{h(a)}A)$ and given an $i\in V(Q)$, define $E_1$ as in Theorem (\ref{resolution}).  We define $p(\nu)_i(v)=\sum_{a\in E_1}\nu_a(av)$.  This sequence is begun by looking at the injective resolution of a quiver algebra without relations found in \cite{twistedquiverbundles}.\medbreak
\begin{cor}
\label{gorenstein}
If $F$ is a consistent tiling on a surface of genus at least one then the associated algebra $A_F$ is Gorenstein in the sense of \cite{dualizingcomplexes}.
\end{cor}
This is especially significant in the genus 1 case, since if $F$ is consistent and $A_F$ is module finite over its centre, and the centre is a Gorenstein Noetherian algebra, we have that for the category of $A_F$-modules, $(-)[3]$ is a Serre functor relative to $\mathcal{K}_{Z(A_F)}$, the dualizing complex for the centre of $A_F$ (Theorem 7.2.14 of \cite{ginz}).  In particular, by viewing finite-dimensional modules over $A_F$ as compactly supported sheaves on $Z(A_F)$ via the forgetful functor, the isomorphism 
\begin{equation}
\label{serreduality}
\mathrm{\textbf{R}Hom}_{Z(A_F)}(\mathrm{\textbf{R}Hom}_{A_F}(M,N),\mathcal{K}_{Z(A_F)})\cong\mathrm{\textbf{R}Hom}_{A_F}(N,M[3])
\end{equation}
yields the \textit{functorial} isomorphism $\mathrm{Ext}^i(M,N)\cong\mathrm{Ext}^{3-i}(N,M)'$ noted in \cite{MR}.  Furthermore, if the centre is a 3 dimensional toric variety, (\ref{serreduality}) restricts to an equivariant functorial isomorphism on the category of $T$-graded modules, where $T$ is the subtorus of $(\Cp^*)^{E(\Gamma)}$ preserving the superpotential relations $\partial/\partial(X_e)W$.  Here the left hand side of (\ref{serreduality}) is taken in the category of $T$-graded modules over $Z(A_F)$, with its $T$-graded dualizing complex.
\section{Dimer configurations}
We expand upon Theorem 5.4 and Corollary 5.5 of \cite{MR}, concerning the link between height functions and finite-dimensional $G_1$-graded cyclic $A_F$-modules, starting with some definitions.
\begin{defn}
A dimer configuration for a graph $\Gamma$ is a subset $D\subset E(\Gamma)$ of the edges of $\Gamma$ such that for every $v\in V(\Gamma)$, a vertex of $\Gamma$, there is a unique $e\in D$ containing $v$.  Two dimer configurations $D', D$ for a graph $\Gamma$ are said to be \textit{asymptotic} if the complements $D \backslash D'$ and $D'\backslash D$ are finite.
\end{defn}
\begin{defn}
Let $F$ be a tiling on the surface $M_g$, $g\geq 1$.  A height function for $\tilde{F}$ is a function from the tiles in $\tilde{F}$ to the integers.  
\end{defn}
\begin{defn}
To an ordered pair $(D,D')$ of asymptotic dimer configurations on the 1-skeleton $\tilde{\Gamma}$ of $\tilde{F}$ we associate the height function $h(D,D')$ in the following way:  the set of $e\in E(\tilde{\Gamma})$ contained in exactly one of $D$ or $D'$ forms a set of disjoint loops $H_{D,D'}$, which we direct so that if an edge in $H_{D,D'}$ comes from $D$ we direct it from black to white, and vice versa for $D'$.  In this way we obtain a set of disjoint \textit{oriented} loops.  We then pick a point $j$ outside all of these loops, and set the height of a tile $i\in\tilde{F}$ to be given by the oriented intersection number of these loops with a path from $j$ to the interior of $i$.
\end{defn}
\begin{defn}
We define an \textit{i-cyclic} $A_F$-module to be a pair $(M,f)$, where $f:e_iA_F\rightarrow M$ is a surjection.  We define a $G_1$-graded $i$-cyclic module to be a pair $(M,f)$ as before, where $M$ is a $G_1$-graded module and $f$ is a surjection of $G_1$-graded $A_F$-modules.  Finally, we define a morphism of $i$-cyclic modules (respectively, $G_1$-graded $i$-cyclic modules) from $(M,f)$ to $(M',f')$ to be a morphism from $f$ to $f'$ in the morphism category of $A_F$ modules (respectively, $G_1$-graded $A_F$-modules).
\end{defn}
\begin{defn}
Let $i$ be a tile of $F$.  After choosing a lift $i_0\in\tilde{F}$ of $i$, a finite-dimensional $i$-cyclic $G_1$-graded $A_F$-module $(M,f)$ gives a height function, denoted $H(M)$, in the following way:   $M$ has a natural decomposition
\begin{equation}
\label{tiledecomp}
M\cong\bigoplus_{j\in \tilde{F}} M_j
\end{equation}
where $M_j$ is spanned as a vector space by $G_1$-homogeneous elements of $M$ whose $G_1$-grade is represented by a path starting at $i_0$ and terminating at $j$.  We define $H(M)(j)=dim(M_j)$.
\end{defn}
Throughout this section we fix, once and for all, a lift $i_0$ of the tile $i$.
\begin{lem}
\label{determine}
Let $F$ be a consistent tiling.  A $G_1$-graded $i$-cyclic $A_F$-module $(M,f)$ is determined up to isomorphism by $H(M)$.
\end{lem}
\begin{proof}
Let $j_0$ be some tile of $\tilde{F}$.  The cyclic module $(M,f)$ is determined up to isomorphism by the kernel of $f$, which we denote $h:I\rightarrow e_i A_F$.  $I$ naturally carries a $G_1$-grading, and $h$ is a $G_1$-graded map.  In particular $I$ admits a decomposition
\[
I=\oplus_{k\in \tilde{F}}I_k.
\]
Let $j_0$ be some $j$-tile of $\tilde{F}$, and consider the summand $I_{j_0}$.  We denote by $\Cp[\omega_j]$ the subalgebra of $A_F$ given, as a vector space, by linear sums of closed $j$-loops.  By consistency, and Lemma (\ref{toruscone}), we have
\[
\mathbb{C}[\omega_j]\cong\mathbb{C}[x].
\]
This algebra carries a natural $G_1$-grading.  Finally, consider the decomposition
\[
e_iA_F=\oplus_{k\in\tilde{F}}e_iA_{F,k},
\]
where $e_iA_{F,k}$ is the direct sum of the $G_1$-graded pieces of $e_iA_F$ corresponding to paths from $i_0$ to $k$ in $\tilde{F}$.  The inclusion $h$ restricts to an inclusion of $G_1$-graded $\Cp[\omega_j]$-modules 
\[
h_{j_0}:I_{j_0}\rightarrow e_iA_{F,j_0}.
\]
Finally, observe that $G_1$-graded submodules of $e_iA_{F,j_0}$ are determined entirely by their colength, since, by consistency, $e_iA_{F,j_0}$ is just a free rank 1 $\Cp[\omega_j]$-module, and this colength is just $H(M)(j_0)$.
\end{proof}
The proof of the following theorem occupies the rest of the section.

\begin{thm}
\label{asymptoticdimers}
For every consistent brane tiling $F$ on a surface of genus at least $1$, and every choice of a tile $i_0\in~\tilde{F}$, a lift of $i\in F$, there is a unique dimer configuration $D_{F,i_0}$ for the 1-skeleton $\tilde{\Gamma}$ of $\tilde{F}$, such that there is a 1-1 correspondence between finite-dimensional $i$-cyclic $G_1$-graded $A_F$-modules $M$ and dimer configurations for $\tilde{\Gamma}$ that are asymptotic to $D_{F,i_0}$.  If $D$ is a dimer configuration for the 1-skeleton $\tilde{\Gamma}$ of $\tilde{F}$ which is asymptotic to $D_{F,i_0}$, the corresponding module is given by taking the unique $i$-cyclic $A_F$-module $M$ satisfying
\[
H(M)=h(D,D_{F,i_0}).
\]
\end{thm}
This is essentially a higher genus version of Theorem 5.4 from \cite{MR}, but as ever we dispense with their second consistency condition and the nondegeneracy condition, since we are using the $G_1$-grading and not the torus grading constructed there.  We will first prove a number of lemmas regarding height functions associated to $i$-cyclic $G_1$-graded $A_F$-modules.\smallbreak
By Lemma \ref{determine}, the isomorphism class of a cyclic module $(M,f)$ is determined by the set
\[
\Omega_M:=\{[u]\hbox{ }|\hbox{ }M_{[u]}\neq 0\},
\]
since this set determines and is determined by the height function of $M$.
\smallbreak
Recall the construction of $h(D_1,D_2)$ for two arbitrary asymptotic dimer configurations: we first obtain from the ordered pair $(D_1,D_2)$ a disjoint union of oriented loops, and then to a disjoint union of oriented loops we associate a height function.  If all the loops are oriented the same way, this height function, applied to some tile in $\tilde{F}$, just counts the minimum number of loops one must cross to get from the interior of the tile to infinity.
\begin{lem}
\label{disjloops}
For all $G_1$-graded $i$-cyclic finite-dimensional $A_F$-modules $M$, $H(M)$ is given by disjoint loops.
\end{lem}
\begin{proof}
First note that this is really a local statement.  If we define $\mathcal{T}_V$ to be the set of tiles containing some vertex $V$ of $\tilde{F}$, then the lemma states that the height function, when restricted to $\mathcal{T}_V$, is either constant, or takes precisely two values $n$ and $n+1$, such that there is a broken line through $\mathcal{T}_V$ consisting of two of the edges of $\tilde{F}$ such that on one side of the broken line $H(M)$ takes value $n$, and on the other side $H(M)$ takes value $n+1$.  If this local statement is true everywhere, the global statement of the lemma easily follows.\smallbreak
We show that $<_{i_0}$ restricts to a total order on $\mathcal{T}_V$ (see Definition (\ref{ordering}) for the definition of $<_{i_0}$).  Let $j$, $k\in\mathcal{T}_V$.  Let $v_{i_0j}$ and $v_{i_0k}$ be minimal paths to the two tiles from $i_0$.  There is a minimal path from $j$ to $k$ consisting of arrows going around the vertex $V$, which we denote $v_{jk}$.  We define $v_{kj}$ similarly.  Clearly we cannot have that $v_{i_0j}v_{jk}$ and $v_{i_0k}v_{kj}$ are both minimal, by antisymmetry of $<_{i_0}$.  So say, without loss of generality, that $v_{i_0k}v_{kj}\sim v_{i_0j}\omega_j^t$, for $t\geq 1$.  We have the further decomposition $v_{i_0j}\omega_j^t\sim v_{i_0j}\omega_j^{t-1}v_{jk}v_{kj}$, from which we deduce, by consistency, that $v_{i_0k}\sim v_{i_0j}\omega_j^{t-1}v_{jk}$.  By minimality of $v_{i_0k}$ we deduce that $t=1$, and that $v_{i_0j}v_{jk}$ is minimal.  We deduce that $j<_{i_0} k$; thus we have proved that $<_{i_0}$ induces a total order on $\mathcal{T}_V$.  \smallbreak
Next note that if $j<_{i_0} k$ then $H(M)(j)\geq H(M)(k)$, since by the assumption $j<_{i_0} k$ we have a minimal path $v_{i_0k}$ with subpath $v_{i_0j}$ a minimal path to $j$, and 
\[\{[v_{i_0k}],[\omega_i v_{i_0k}],\ldots,[\omega_i^{H(M)(k)-1}v_{i_0k}]\}\subset \Omega_M,
\] and so 
\[\{[v_{i_0j}],[\omega_i v_{i_0j}],\ldots,[\omega_i^{H(M)(k)-1}v_{i_0j}]\}\subset \Omega_M.
\]
Similarly, since $j$ and $k$ share a vertex $V$, the difference between $H(M)(j)$ and $H(M)(k)$ is at most 1, since if $H(M)(j)=t$ then 
\[
\{[v_{i_0j}],[v_{i_0j}\omega_j],\ldots,[v_{i_0j}\omega_j^{t-1}]\}\subset\Omega_M
\]
and so 
\[\{[v_{i_0j}v_{jk}],[v_{i_0j}\omega_j v_{jk}],\ldots,[v_{i_0j}\omega_j^{t-2}v_{jk}]\}\subset\Omega_M.
\]  
\smallbreak
We have shown that $<_{i_0}$ induces a total order of $\mathcal{T}_V$.  Let $k_1$ be the maximum with respect to this ordering.  Let $k_2$ be the maximum element of $\mathcal{T}_V\backslash k_1$, with respect to $<_{i_0}$.  Then there is a minimal path $u$ from $i$ to $k_1$ such that $u$ decomposes as $u=vw$, where $v$ is a minimal path from $i_0$ to $k_2$ and $w$ is a minimal path from $k_2$ to $k_1$.  Since $k_1$ and $k_2$ share a vertex $V$, there is a minimal path from $k_2$ to $k_1$ only crossing those edges of $\tilde{F}$ that contain $V$.  It follows that the path $w$ consists of a solitary arrow, for else there is some $k_3\neq k_1$ such that $k_3$ contains $V$, and such that $k_2<_{i_0} k_3$.  We deduce that, assuming the dual quiver $\tilde{Q}$ is oriented clockwise around $V$ (i.e. $V$ is coloured black), the tiles are arranged, in order with respect to $<_{i_0}$, clockwise.  If, on the other hand, $V$ is coloured white, then the tiles are arranged, in order with respect to $<_{i_0}$, anticlockwise.
\smallbreak
The function $H(M)$ is monotone decreasing with respect to the ordering $<_{i_0}$, and the maximum difference between $H(M)(j)$ and $H(M)(k)$ for $k,j\in\mathcal{T}_V$ is 1.  It follows that, at $V$, the height function is given by a pair of lines through $V$, passing along two of the edges of $\tilde{F}$ containing $V$, or by no lines at all.  This implies that the height function is given by a disjoint union of loops.
\end{proof}
\smallbreak
Now, let $H(M)(i_0)=n$, for some $G_1$-graded $i$-cyclic module $(M,f)$.  Then consider the $G_1$-graded short exact sequence
\begin{equation}
\label{norespect}
\xymatrix{
0\ar[r]&Image(f|_{[\omega_i]^{n-1}A_F})\ar[r]&M\ar[r]&N\ar[r]&0.
}
\end{equation}
There is a natural cyclic module $(Image(f|_{[\omega_i]^{n-1}A_F}),g)$ with $g$ sending $1\in A_F$ to $f([\omega_i]^{n-1})$.  The naturally defined height function $H(Image(f|_{[\omega_i]^{n-1}A_F}))$ has maximum value $1$, since it takes value $1$ at $i_0$.  If an $i$-cyclic module $(M,f)$, with underlying module $M$ possessing a $G_1$-grading, satisfies
\[
\mathrm{image}(H(M))\subset[0,1]
\] 
we call $M$ a \textit{bungalow}.  An $i$-cyclic submodule $(M',f')$ of $(M,f)$, with $M'$ a $G_1$-graded submodule of $M$ generated in a grade lying above the tile $i_0$, is determined up to isomorphism by the minimum value of $r$ such that $[\omega_i]^r\in M'$ (identifying $M$ with a quotient of $e_i A_F$).  If $(M',f')$ is a bungalow this value is forced, and must be $n-1$.  We deduce that $(Image(f|_{[\omega_i]^{n-1}A_F}),g)$ is the unique (up to isomorphism) submodule of $M$ that is a bungalow.  If 
\[
M'\rightarrow M\rightarrow M''
\]
is a $G_1$-graded short exact sequence of underlying modules of cyclic modules, it follows from the construction of the height functions associated to the three modules that $H(M)=H(M')+H(M'')$.  Note any quotient of the underlying module of a cyclic module naturally has the structure of a cyclic module.  We deduce, by induction on $H(M)$, that a $G_1$-graded finite-dimensional cyclic module has a unique filtration by bungalows generated in grade lying above tile $i_0$.  We will need the following lemma:
\begin{figure}
\centering
\input{newdig.tex}
\caption{winding around the loop $\mathcal{L}$}
\label{newdig}
\end{figure}
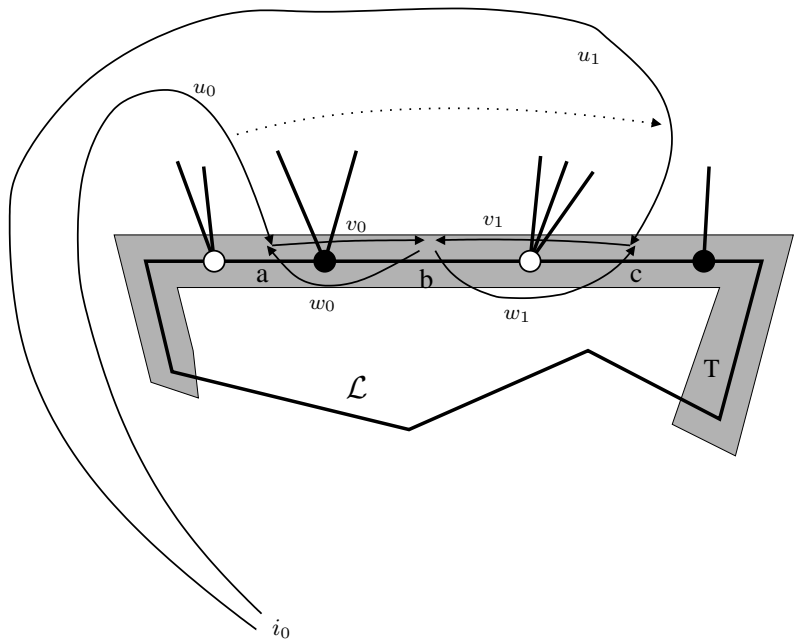

\begin{lem}
The height function $H(M)$ associated to a bungalow $(M,f)$ is given by a single loop.
\end{lem}
\begin{proof}
Let $H_1$ be the set of tiles in $\tilde{F}$ sent to $1$ by $H(M)$.  $H_1$ is connected, since $H(M)$ is monotone decreasing with respect to $<_{i_0}$.  We need to show that it is simply connected.  Unsurprisingly, this will involve winding.  Let $\mathcal{L}$ be a loop of edges of $\tilde{\Gamma}$ bounding a hole in $H_1$ - i.e. the tiles on the outside of $\mathcal{L}$ containing edges of $\mathcal{L}$ are in $H_1$, and the tiles inside $\mathcal{L}$ are not contained in $H_1$.  Since $H(M)$ is the height function associated to a module, it follows by Lemma (\ref{disjloops}) that it is given by disjoint loops, and so all of the tiles outside $\mathcal{L}$ with boundary intersecting $\mathcal{L}$ nontrivially are contained in $H_1$.  Let $a,b,c$ be three arrows of $\tilde{Q}$ as indicated in Figure (\ref{newdig}) (we label the edges of $\tilde{F}$ dual to these arrows).  Let $u_0$ be a minimal path from $i_0$ to $h(a)$.  There is a  minimal path $v_0$ from $h(a)$ to $t(b)$ as indicated in Figure (\ref{newdig}).  We claim that $u_0v_0$ is minimal.  To see this, assume otherwise.  Then $u_0v_0\sim w\omega_{t(b)}$ for some path $w$.  We can pick $\omega_{t(b)}$ so that it decomposes as $w_0v_0$, where $w_0$ is as in Figure (\ref{newdig}).  By consistency it follows that $u_0\sim w w_0$.  since $H(M)$ is monotone decreasing with respect to $<_{i_0}$ we deduce that $w w_0$ passes only through tiles contained in $H_1$, but we picked $w_0$ to leave $H_1$, a contradiction.  We deduce that $u_0v_0$ is minimal.\smallbreak
Let $u_1$ be a minimal path to $h(c)$.  By exactly the same argument, if $v_1$ is the minimal path from $h(c)$ to $t(b)$ indicated in Figure (\ref{newdig}) then $u_1v_1$ is minimal.  By consistency, all minimal paths are equivalent, and so there is a chain of basic F-term equivalent paths taking us from $u_0v_0$ to $u_1v_1$.  These paths must all be contained in $H_1$, since $H(M)$ is monotone decreasing with respect to $<_{i_0}$.  
\smallbreak
Recall the definition of $\Psi$ from (\ref{psidef}) of Section \ref{homotopies}, a function from suitably nice continuous paths in $M_g$ to paths in $\tilde{F}$.  Let $T$ be a tubular neighborhood of $\mathcal{L}$ (part of $T$ is indicated in Figure (\ref{newdig})).  Let $\gamma_i$ satisfy $\Psi(\gamma_i)=u_i$ for $i=0,1$.  Let $\lambda_i$ satisfy $\Psi(\lambda_i)=v_i$ for $i=0,1$.  As indicated, we can choose $\lambda_i$ to lie entirely within $T$, and for $\gamma_i$ to end at the starting point of $\lambda_i$, for each $i$.  Then we can construct a homotopy from $\gamma_0$ to $\gamma_1$ such that $image(H|_{[0,1]\times\{1\}})\subset T$, by the following procedure.  First, take a homotopy that stretches $\gamma_0$ to the path given by concatenating $\gamma_0$ and $\lambda_0$.  Similarly, there is a homotopy that stretches $\gamma_1$ along the path $\lambda_1$ marked in Figure (\ref{newdig}).  Finally, by F-term equivalence of $u_0v_0$ and $u_1v_1$, there is an endpoint preserving homotopy between these two stretched paths, which avoids $\mathcal{L}$. We can pick $T$ to be arbitrarily thin.  We have shown that there is a $\gamma_0$ such that $\Psi(\gamma_0)=u_0$, and a $\gamma_1$ such that $\Psi(\gamma_1)=u_1$, and a homotopy
\[
H:[0,1]^2\rightarrow H_1
\]
such that $H|_{\{i\}\times[0,1]}=\gamma_i$, and $im(H|_{[0,1]\times\{1\}})\subset T$.  After repeatedly applying this argument we arrive at a $\gamma_{\infty}$ such that $\Psi(\gamma_{\infty})\sim u_0$, but $\gamma_{\infty}$ doesn't have the same winding number around $\mathcal{L}$ as $\gamma_0$, a contradiction since $\gamma_0,\gamma_{\infty}$ are $k$-unwound for any tile $k$ not in $H_1$.  This proves the lemma.
\end{proof}

So far we have proved that an $i$-cyclic $G_1$-graded module $M$ is determined by its filtration by bungalows.  Furthermore, a bungalow determines a single loop in $\tilde{F}$, defined by its height function.  The height function for $M$ is given by the disjoint union of these loops, and $M$ is uniquely determined by this disjoint union of loops.  In one direction, we have shown that a module determines a finite disjoint union of loops that satisfy the property that they are the height functions of bungalows.  By the next lemma, this is a 1-1 correspondence
\begin{lem}
If $H$ is the height function associated to a finite disjoint union of loops arising from bungalows, then $H$ is the height function associated to a $G_1$-graded $i$-cyclic module.
\end{lem}
\begin{proof}
Let $\mathcal{L}_s$, $s\in [1,\ldots,n]$ be a finite union of disjoint loops, each the height function of a bungalow $M_s$.  Since they are disjoint, and all contain $i_0$, they form a series of concentric loops around $i_0$.  Define 
\[
\Omega:=\{[v_{i_0j}][\omega_j]^r|\text{ for }j\text{ inside at least }r-1\text{ of the loops of }\mathcal{L}\}.
\]
Since each loop $\mathcal{L}_s$ of $\mathcal{L}$ is the height function of a module, it follows that if $j_0$ is a tile contained in $\mathcal{L}_s$ and $u$ is a minimal path from $i_0$ to $j_0$, $u$ is contained entirely in $\mathcal{L}_s$.  We deduce that if $u$ is a path starting at $i_0$ that crosses $n$ of the loops at least twice, then $u$ factors as $u\sim v\omega_{h(u)}^n$ for some minimal $v$.  It follows that the set of equivalence classes of paths that are not in $\Omega$ is closed under extension of paths, and so $\Omega$ is the basis for an $i$-cyclic module.
\end{proof}
Let $D$ be a dimer configuration on $\tilde{\Gamma}$, the 1-skeleton of $\tilde{F}$.  We say that an oriented loop $\mathcal{L}$ is \textit{obtainable} from $D$ if alternating edges of $\mathcal{L}$ belong to $D$, and the edges belonging to $D$ are directed from black to white vertices.  The construction of $h(D,-)$ gives a correspondence between dimer configurations that are asymptotic to $D$, and finite sets of obtainable loops of $D$.  We are finally in a position to prove Theorem (\ref{asymptoticdimers}), since we have reduced it to the following proposition:
\begin{prop}
For every consistent brane tiling $F$ on a surface of genus at least $1$, and every choice of a tile $i_0\in~\tilde{F}$, there is a unique dimer configuration $D_{F,i_0}$ for the 1-skeleton of $\tilde{F}$, such that the obtainable loops of $D_{F,i_0}$ are exactly the loops corresponding to height functions of bungalows.
\end{prop}

\begin{figure}[htbp]
\begin{center}
\caption{The bold edges show those from $D_n$, as well as new ones for $D_{n+1}$.}
\label{loopextension} 
\input{expansion.tex}
 
\end{center}
\end{figure}
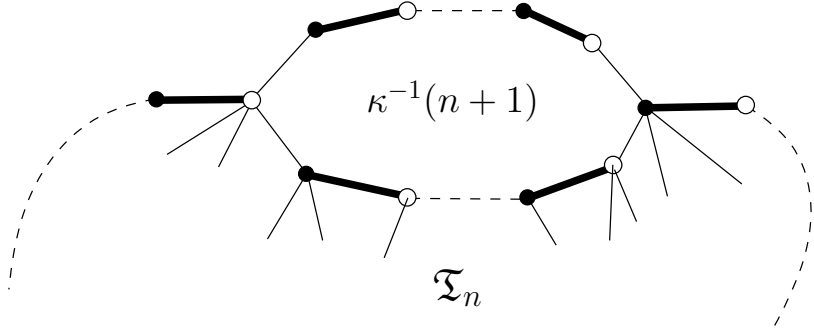

\begin{proof}
We extend the partial order $<_{i_0}$ to a total order, which we will denote $<$, and obtain in this way a function $\kappa$ from the tiles in $\tilde{F}$ to $\mathbb{N}$.  Let $\mathfrak{T_n}=\kappa^{-1}([1,\ldots,n])$.  We assume that we have constructed a dimer configuration $D_n$ on $\mathfrak{T_n}$ satsifying the condition above - the set of obtainable loops of $D_n$ in the 1-skeleton of $\mathfrak{T_n}$ is in bijective correspondence with the set of loops in the 1-skeleton of $\mathfrak{T_n}$ that arise from bungalows.  For $n=1$, existence and uniqueness of this configuration is clear.\smallbreak
Note that the boundary of $\mathfrak{T_n}$ is the loop corresponding to the bungalow which is the cokernel of the map
\[
\xymatrix{
[\omega_i] A_F \oplus \bigoplus_{j\notin \mathfrak{T_n}}[v_{i_0j}] A_F\ar[r]&e_iA_F
}
\]
and so the correct alternating edges of the boundary of $\mathfrak{T}_n$ belong to $D_n$, by induction.  Figure (\ref{loopextension}) depicts a part of the loop bounding $\mathfrak{T}_n$, along with the tile $\kappa^{-1}(n+1)$.  Recall from the proof of Lemma (\ref{disjloops}) that for each vertex $j$, the restriction of $<_{i_0}$ to the tiles containing $j$ is already a total order, and around the white vertices the tiles are arranged anticlockwise, with respect to this ordering, and around the black vertices they are arranged clockwise.  It follows that the addition of the new tile must take the form as depicted in Figure (\ref{loopextension}), with the part of the loop bounding $\mathfrak{T}_{n+1}$ that is not contained in $\mathfrak{T}_n$ being a path which, travelling clockwise, goes from a white to a black vertex. Since we require the boundary of $\mathfrak{T}_{n+1}$ to be a loop obtainable from $D_{n+1}$, from Figure (\ref{loopextension}) it is clear how we \textit{must} (hence, uniqueness) extend $D_n$, along the edges marked in bold, where the edges that must already be in $D_n$ are also marked in bold.  
\smallbreak
We now give the inductive step: assume we have a loop $S$, obtainable from $D_{n+1}$, contained in the 1-skeleton of $\mathfrak{T}_{n+1}$.  We assume that it encloses $\kappa^{-1}(n+1)$.  Then from the figure if we modify the loop by removing $\kappa^{-1}(n+1)$ from its interior, we obtain a loop in the 1-skeleton of $\mathfrak{T}_n$ obtainable from $D_n$.  By induction this loop $S'$ corresponds to a bungalow given by the cokernel of the map
\[
\xymatrix{
[\omega_i] \oplus \bigoplus_{j\notin Int(S')} [v_{i_0j}]A \ar[r] &e_i A.
}
\]
Any minimal path to $\kappa^{-1}(n+1)$ must pass through one of the tiles in $\mathfrak{T}_n$ adjoining $\kappa^{-1}(n+1)$, and so must lie entirely within $S'$, until it leaves for $\kappa^{-1}(n+1)$.  It follows that the set $P$ of paths $u\in paths_{\tilde{F}}$ that are minimal paths from $i_0$ to tiles in the interior of $S$ is closed under the operation $u\rightarrow u'$ of deleting the final arrow, and so if $Q$ is the set of F-term equivalence classes of paths not in $P$, $Q$ spans a right ideal of $A_F$, which we call $I_Q$.  By construction, the cokernel of
\[
\xymatrix{
I_Q\ar[r]&e_i A
}
\]
gives the required module.  Conversely, say we have a bungalow $M$ with associated loop $S$ contained in $\mathfrak{T}_{n+1}$.  Without loss of generality it encloses $\kappa^{-1}(n+1)$, and we set $j=\kappa^{-1}(n+1)$.   Then the bungalow given by the cokernel of the map
\[
\xymatrix{
[v_{i_0j}]A\ar[r]&M
}
\]
gives a loop $S'$ contained in $\mathfrak{T}_n$, and, by induction, alternating edges are contained in $D_n$.  Furthermore, since $<_{i_0}$ is a total order when restricted to the set of tiles containing any given vertex of $\tilde{F}$, it follows that $S'$ contains all the tiles of $\mathfrak{T}_n$ that border $\kappa^{-1}(n+1)$.  We deduce from the indicated extension of $D_n$ that the loop $S$ has alternating edges in $D_{n+1}$.  This completes the proof of the proposition, and the proof of Theorem (\ref{asymptoticdimers}).
\end{proof}

\section{Donaldson-Thomas theory}
The final section of \cite{MR} deals with the Donaldson-Thomas theory of the moduli space $\mathcal{M}_{\textbf{v},i}$ constructed in \cite{conifold}, of $i$-cyclic modules with dimension vector $\textbf{v}\in \mathbb{N}^{V(Q)}$.  Given an action of a torus $T$ on the moduli space of cyclic modules, such that the perfect obstruction theory on the moduli space is equivariant with respect to $T$, and the fixed locus consists of isolated fixed points, Corollary 3.5 of \cite{ObstructionsHilbert} states that the Donaldson-Thomas invariant of the moduli space is given by the sum, over the $T$-fixed points, of the sign of the dimension of the Zariski tangent space at those points.  In order to apply Corollary 3.5 of \cite{ObstructionsHilbert} we need to find the fixed points of the moduli space, under some action of a torus $T$ that respects the perfect obstruction theory, confirm that the Zariski tangent space has no trivial weights under the induced action of $T$, and then compute the dimensions of the Zariski tangent spaces.\smallbreak
If $T'\subset (\Cp^*)^{|E(F)|}$ is the subtorus preserving the superpotential relations, then there is a natural action of $T'$ on $\mathcal{M}_{\textbf{v},i}$.  This action is not faithful, and the kernel of the map 
\begin{equation}
\label{theta}
\theta:T'\rightarrow \prod_{\textbf{v}\in\mathbb{N}^{V(Q)}}\mathrm{Aut}(\mathcal{M}_{\textbf{v},i})
\end{equation}
 is given by the subtorus generated by those 1-parameter subgroups of $T'$ that rescale all arrows going into one vertex by some number $z\in\Cp$, and rescale all arrows going out by $z^{-1}$.  We let $T$ be the image of $\theta$ (again $T$ is a torus).  Let $T_W$ be the image under $\theta$ of the largest subtorus of $T'$ acting with trivial weight on the simple loop $\omega_i$.  The obstruction theory constructed in \cite{conifold} is symmetric equivariant with respect to the induced action of $T_W$.  
\begin{lem}
The torus $T_W$ has dimension $2g$.
\end{lem}
\begin{proof}
We show the coweight lattice of $T_W$ is $\mathbb{Z}^{2g}$.  Pick $i$ a tile of $F$, and let $i_0\in\tilde{F}$ be an $i$-tile.  Let 
\[
\alpha:\Cp^*\rightarrow T_W
\]
be a coweight.  Then $\alpha$ is determined by a $\Cp^*$-action on $A_F$ that acts trivially on $\omega_i$.  Note that by transport of simple loops (see Remark (\ref{transport})) this implies $\Cp^*$ acts trivially on all simple loops.  Such an action is determined by its effect on minimal paths.  Since we have picked an initial $i_0$, it makes sense to talk of the action of $\Cp^*$ on a tile in $\tilde{F}$ under $\alpha$, and in fact $\alpha$ is determined by the actions of $\Cp^*$ on the tiles of $\tilde{F}$.  For any $u$ a $j$-loop, there is a $v$ which is an $i$-loop, and a path $w$ from $i$ to $j$, such that $wu\omega_j^n\sim vw$ for some $n\geq 0$.  It follows that the $\Cp^*$-action on $j$-loops is determined by the $\Cp^*$-action on $i$-loops.  Furthermore, increasing the weights on all $j$-tiles  corresponds to an action of the gauge group.  It follows that the action of $\Cp^*$ on $\prod_{\textbf{v}\in\mathbb{N}^{V(Q)}}\mathrm{Aut}(\mathcal{M}_{\textbf{v},i})$ determines and is determined by its action on the $i$-tiles of $\tilde{F}$.  Such an action uniquely determines an element of $Hom(H_1(M_g,\mathbb{Z}),\mathbb{Z})$, since for any two $i$-loops $v$ and $w$, the action on the endpoint of the lift $\widetilde{vw}$ starting at $i_0$ is the same as the action on the endpoint of the lift $\widetilde{wv}$.  So all that is required is to show that every element of $Hom(H_1(M_g,\mathbb{Z}),\mathbb{Z})$ is realised in this way.
\smallbreak
Assume that we have drawn the dual graph of $F$ on $M_g$, so we have some embedding $Q\subset M_g$.  Let $C$ be a curve in $M_g$ representing some $\beta\in H_1(M_g,\mathbb{Z})$.  We can pick $C$ so that it doesn't go through the vertices of $Q$.  Let $p^{-1}(C)$ be the preimage of $C$ under the projection from the universal cover.  Now number the tiles in $\tilde{F}$ so that the difference between two adjacent tiles $j$ and $k$ is given by the oriented intersection number of the arrow between the vertices in the dual graph corresponding to $j$ and $k$, with $p^{-1}(C)$.  We give tile $i_0$ the number 0.  This gives an assignment of numbers to the tiles of $\tilde{F}$.  By construction, the corresponding element of $Hom(H_1(M_g),\mathbb{Z})$ is the image of $\beta$ under the natural isomorphism induced by the nondegenerate intersection product.  This completes the proof.
\end{proof}

\begin{prop} If $F$ is a consistent tiling on a surface of genus 1, then for all $\textbf{v}\in\mathbb{N}^{V(Q)}$ the $T_W$-fixed points of $M_{\textbf{v},i}$ are the $G_1$-graded quotients of $e_iA$ with dimension vector $\textbf{v}$.
\label{fixedpoints}
\end{prop}
\begin{rem}
In fact the above holds if and only if $F$ is consistent.
\end{rem}
\begin{proof}
Let $(M,f)$ be an $i$-cyclic $A_F$-module, and let 
\[
I\rightarrow A_F\rightarrow M
\]
be the associated short exact sequence.  We consider $I$ as a right ideal of $A_F$.  $I$ is generated by elements of the form 
\begin{equation}
\label{Nform}
\sum_{s\in S} \lambda_s [u_s],
\end{equation}
where the $\lambda_s\in \Cp$ and $u_s$ are paths from $i$ to some fixed $j\in F$.  Since $M$ is $T_W$ invariant, $I$ is generated by elements of the form in (\ref{Nform}), chosen such that the lifts of the $u_s$ to paths starting at $i_0$ in $\tilde{F}$ all terminate at the same lift $j_0$ of $j$.  We deduce that, as a vector space, we have a decomposition
\begin{equation}
\label{tiledecomp2}
I\cong\bigoplus_{k_0\in \tilde{F}} I_{k_0}
\end{equation}
where each $I_{k_0}$ is a sub-$\Cp[\omega_{k}](\cong\Cp[x])$-module of $e_iA_{F,k_0}$, where $e_iA_{F,k_0}$ is the summand associated to the tile $k_0$ in the analogous decomposition of $e_iA_F$.  By consistency, each $e_iA_{F,k_0}$ is isomorphic to the free rank 1 $\Cp[x]$-module.  All that is required is to show that in fact each $I_{k_0}$ is a $\mathbb{N}$-graded submodule of $e_iA_{F,k_0}$.\smallbreak
We proceed by induction on the codimension of $I_{k_0}$ in $e_iA_{F,k_0}$.  If it is zero, there is nothing to prove.  Assume it is $n$, and we have proved the result for all $k_0$ such that the codimension of $I_{k_0}$ is less than $n$.  Let $u$ be a minimal path from $i_0$ to $k_0$.  By Lemma (\ref{MRres}) of the next section, we deduce that there is an arrow $a_0$ such that $ua_0$ is minimal.  Denote $h(a_0)$ by $k_1$.  The ideal $I_{k_0}\subset e_iA_{F,k_0}$ is generated by some degree $n$ polynomial $p(x)$ (in $\omega_{k_0}$), i.e. $up(\omega_{k_0})\in I$.   Since $I$ is a right ideal, we deduce that $ua_0p(\omega_{k_1})\in I$, and the same polynomial belongs to $I_{k_1}$.  We deduce that the codimension of $I_{k_1}$ is less than or equal to $n$.  If it is equal to $n$, then in fact right multiplication by $a_0$ induces an isomorphism from $I_{k_0}$ to $I_{k_1}$, since any other degree $n$ polynomial contained in $I_{k_1}$ must be a multiple of $p$.  We obtain in this way a chain of isomorphisms until, by finiteness and the fact that minimal paths do not revisit tiles, we arrive at some $k_t$, and arrow $a_t$ from $k_t$ to some $k_{t+1}$, such that the codimension of $I_{k_t}$ is strictly greater than the codimension of $I_{k_{t+1}}$.  $I_{k_{t+1}}$ is an ideal in $e_iA_{F,k_{t+1}}$, generated by some degree $l$ polynomial $q$.  By induction, this polynomial is $q(x)=x^{l}$.  Right multiplication by $\omega^{br}_{a_{t}}$ defines a $\Cp[x]$-module map from $I_{k_{t+1}}$ to $I_{k_{t}}$, which takes $ua_0\ldots a_t\omega_{k_{t+1}}^l$ to $ua_0\ldots a_{t-1}\omega_{k_t}^{l+1}$.  Since we choose $ua_0\ldots a_{t-1}$ to be minimal, this represents the element $x^{l+1}$ under the isomorphism $e_iA_{F,k_t}\cong\Cp[x]$.  Since $I_{k_t}$ has codimension $n$ in $e_iA_{F,k_t}$, it follows that $n\leq l+1$ and so $l=n-1$, since $l<n$.  It follows that $p(x)$ is a scalar multiple of $x^n$.
\end{proof}
We next prove the following:
\begin{prop}
\label{trivialweights}
If $F$ is a tiling on a surface of genus 1, and $M_{\textbf{v},i}$ is the moduli space associated to the dimension vector $\textbf{v}$ then the action of $T_W$ on the Zariski tangent space at the fixed points carries no trivial representations.
\end{prop}

\begin{proof}
Let $I_0$ be the kernel of the map $f:e_i A\rightarrow M$.  We consider $I_0$ as a right ideal of $A_F$.  As in Lemma 4.1(b) of \cite{ObstructionsHilbert}, we check the condition via the identification of the Zariski tangent space with the space $\mathrm{Hom}(I_0,M)$.  Note that both $I_0$ and $M$ are $G_1$-graded by Proposition (\ref{fixedpoints}).\smallbreak
Let $T'_W\subset T'$ be the subtorus of $(\Cp^*)^{|E(Q)|}$ preserving the superpotential relations and acting with trivial weight on the simple loops.  $\mathrm{Hom}(I_0,M)$ carries a natural $T'_W$-action, given by the action of $T'_W$ on the tangent space of a fixed point of the $T'$-action on $M_{\textbf{v},i}$.  This coincides with the natural $T'_W$-action given on homomorphisms between modules equipped with a $T'_W$-action.  The $T'_W$-trivial sub-representations are just the $T_W$-trivial sub-representations, by the definition of $T_W$. 
\smallbreak
Any $A_F$-module homomorphism from $I_0$ to $M$ takes paths ending at $j$-tiles to paths ending at $j$-tiles.  A homomorphism is acted on trivially by $T_W$ (and hence by $T'_W$) if and only if it takes elements of the form $[v_{i_0j_0}\omega_{j_0}^n]$, where $v_{i_0j_0}$ is some minimal path from $i_0$ to $j_0$ in $\tilde{F}$, to linear combinations of elements of the form $[v_{i_0j_0}\omega_{j_0}^m]$, with the same $j_0\in\tilde{F}$.  Now let $\phi\in\mathrm{Hom}(I_0,M)$ be a nonzero homomorphism acted on trivially by $T_W$.\smallbreak
Recall from the proof of Theorem (\ref{asymptoticdimers}) that the module $M$ is equivalent to the data of its height function, and this is given by a set of disjoint loops, each one the height function of one of the modules occuring in its filtration by `bungalows'.  These loops form a concentric series of loops around $i_0$, since they are disjoint and all contain $i_0$.  Let the loops be ordered so that the outermost is $\mathcal{L}_1$, i.e. $\mathcal{L}_t$ encloses exactly those tiles $j_0$ such that $H(M)(j_0)\geq t$.  Let $j_0$ lie outside $\mathcal{L}_1$.  Since $\phi$ is acted on trivially by $T_W$ it preserves the endpoints of paths.  It follows that $\phi([v_{i_0j_0}])=0$, where $v_{i_0j_0}$ is the minimal path from $i_0$ to $j_0$, since the $G_1$-graded pieces of $M$ corresponding to paths ending at $j_0$ are all zero.  Similarly, for any $[v]\in I_0$ with $h(v)=j_0$, $\phi([v])=0$.  
\smallbreak
Say, by induction, that for all $j_0$ outside $\mathcal{L}_t$ it is shown that, for all $[v]\in I_0$ with $h(v)=j_0$, $\phi([v])=0$.  Let $j_0$ be a tile on the interior of $\mathcal{L}_t$, sharing an edge with $\mathcal{L}_t$.  Let $k_0$ be the other tile sharing this edge (note that $k_0$ lies outside $\mathcal{L}_t$, and $H(M)(k_0)=t-1$).  Then the two minimal paths $v_{k_0j_0}$ and $v_{j_0k_0}$ satisfy $v_{j_0k_0}v_{k_0j_0}\sim \omega_{j_0}$.  By the definition of $\mathcal{L}_{t}$, $[v_{i_0j_0}][\omega_{j_0}^{t-1}][v_{j_0k_0}]\in I_0$, and by induction $\phi([v_{i_0j_0}][\omega_{j_0}^{t-1}][v_{j_0k_0}])=0$.  Any $[v]\in I_0$ with $h(v)=j_0$ can be written $[v_{i_0j_0}][\omega_{j_0}^{t-1}][v_{j_0k_0}][v_{k_0j_0}][\omega_{j_0}^r]$, for some $r\geq 0$.  It follows that for all $[v]\in I_0$ with $h(v)=j_0$, $\phi([v])=0$.  
\smallbreak
Now let $j_0$ and $k_0$ be two arbitrary neighbouring tiles, both lying inside $\mathcal{L}_t$ and outside $\mathcal{L}_{t+1}$.  Let $a$ be the arrow between them, without loss of generality it goes from $j_0$ to $k_0$.  Either $v_{i_0j_0}a$ is minimal, or $v_{i_0k_0}\omega^{br}_{a}$ is.  Consider the decomposition
\[
M\cong \bigoplus_{l\in \tilde{F}}M_l
\]
where $M_l$ is spanned by those $G_1$-graded pieces of $M$ corresponding to paths that end at $l$.  If $v_{i_0j_0}a$ is minimal, then right multiplication by $a$ defines an isomorphism from $M_{j_0}$ to $M_{k_0}$.  Similarly, if $v_{i_0k_0}\omega^{br}_{a}$ is minimal then right multiplication by $\omega^{br}_{a}$ defines an isomorphism from $M_{k_0}$ to $M_{j_0}$.  Either way, we deduce that $\phi([v_{i_0j_0}\omega_{j_0}^{t+1}])=0$ if and only if $\phi([v_{i_0k_0}\omega_{k_0}^{t+1}])=0$.  This completes the inductive step, and the proposition follows.
\end{proof}
\smallbreak
Given these two propositions, we can compute Donaldson-Thomas invariants, assuming only consistency, thanks to the characterisation of the $G_1$-graded cyclic modules of Theorem (\ref{asymptoticdimers}).  We have proved all the essential results of \cite{MR} assuming only consistency.

\section{Other notions of consistency}
We recall the definitions found in \cite{MR}.
\begin{defn}
\label{cons2}
A brane tiling $F$ satisfies `MR2' if, for any two tiles $i,j\in \tilde{F}$, there exists an arrow $a$ from $j$ to some tile $k$ such that some minimal path from $i$ to $k$ passes along $a$.  We also require the dual to hold.
\end{defn}
\begin{defn}
\label{nondegeneracy}
A brane tiling $F$ is \textit{nondegenerate} if its 1-skeleton $\Gamma$ satisfies the following condition: for each edge $E\in E(\Gamma)$ there is a dimer configuration for $\Gamma$ containing $E$.
\end{defn}
The consistency condition of Definition (\ref{consistency}), and the two conditions above, are assumed in \cite{MR}, and are used to prove the results there.  It is natural to wonder whether they follow from the consistency condition (the `first consistency condition' in \cite{MR}).  The answer, at least for the second consistency condition, is affirmative.
\begin{lem}
\label{MRres}
Consistency implies condition MR2.
\end{lem}
\begin{proof}[Sketch proof:] Pick $i$, $j\in \tilde{F}$ violating the condition.  Define $E_1$ as the set of arrows $a$ with $h(a)=j$, and $E_2$ as the set of arrows with $t(a)=j$.  For each $k,l\in \tilde{F}$ let $v_{kl}$ denote a minimal path from $k$ to $l$.  Then we are assuming that for each $a\in E_2$, $v_{ij}a\sim v_{ih(a)}\omega_{h(a)}^t$ for some $t\geq 1$.  If $t>1$ we pick $\omega_{h(a)}$ going through $j$ and deduce that $v_{ij}$ is not minimal.  It follows that for all $a\in E_2$, $v_{ij}\sim v_{i,h(a)}\omega^{br}_{a,B}\sim v_{ih(a)}\omega^{br}_{a,W}$.  None of these paths go through $j$ (except to end at it) by minimality of $v_{ij}$.  We are now in exactly the situation of the winding argument of (\ref{resolution}), giving a contradiction.  The dual argument is identical.
\end{proof}
So condition MR2 is redundant, and if we include our consistency condition we can drop it from our definition of a `good' algebra.  While it is possible that nondegeneracy follows from consistency too, we have no proof of this (though see \cite{kazushi}).\smallbreak

\section{Final remarks}
\begin{thm}
Let $F$ be a consistent brane tiling on an oriented surface $M_g$.  The associated algebra is 3-Calabi-Yau if and only if $g>0$.
\end{thm}
\begin{proof}
From Proposition (\ref{exactCY}), and Theorem (\ref{resolution}), it is enough to show that the complex (\ref{res1}) is not a resolution of $S_i$ in the genus zero case.  For this, it is sufficient to show that one of the $G_1$-graded pieces corresponding to a path of length at least 1 is not exact.  Let $F$ be a consistent brane tiling on $S^2$.  From the proof of Theorem (\ref{resolution}) it is enough to show, using the same notation, that there is a pair of tiles $i$ and $j$ such that for $u$ a minimal path between them $E_1'=E_1$, i.e. that we are in situation (B) of that proof.\smallbreak
In the case of the sphere we have $F=\tilde{F}$, and Lemma (\ref{toruscone}) applies to $F$.  Now for two tiles $i$ and $j$ of $F$ let $v_{ij}$ be a minimal path between them, and let $S(i,j)$ be the set of tiles contained in paths that are F-term equivalent to $v_{ij}$.  Say we have picked $i$ and $j$ so that $S(i,j)$ is maximal.  Note that, for $a\in E_1$, $t(a)$ is in $S(i,j)$ if and only if there is a minimal $v\sim v_{ij}$ starting with $\omega^{br}_{a}$, i.e. for $a\in E_1$, $a\in S(i,j)$ if and only if $a\in E_1'$.  Say there is some $a\notin E_1'$.  It follows that $av_{ij}$ is minimal, since if $av_{ij}\sim\omega_{t(a)}v\sim a \omega^{br}_{a}v$ for some $v$, we deduce that $a\in E_1'$.  But now clearly $|S(t(a),j)|>|S(i,j)|$, a contradiction.
\end{proof}
\begin{rem}
\label{sphere}
Another topological argument shows that in fact $S(i,j)$ contains all the tiles of $F$.  Briefly, we use the same argument as above to deduce that all the tiles sharing edges with $j$ are contained in $S(i,j)$.  Now, say $k\notin S(i,j)$ for some $k\in F$.  For each tile bordering $j$ there is a tile bordering $i$, with some $w\sim \:'v_{ij}'$ going between the two, for some choice of minimal $v_{ij}$.  We wind around $j$ using these paths, necessarily winding around $i$ too as we go ($S^2\backslash\{i,j\}$ is an annulus).  After winding round $j$ we obtain two F-term equivalent paths from a tile bordering $j$ to a tile bordering $i$, which form a cycle around the hole left by $k$ in $S(i,j)\backslash\{i,j\}$.  This proof is suggestive: a consistent tiling on a surface $M_g$ is not Calabi-Yau if and only if we can find the topological obstruction to $M_g$ being a $K(\pi,1)$ surface (i.e. $\pi_2\neq 0$), in one of the $S(i,j)$.
\end{rem}
\begin{rem}
Let $F$ be a consistent brane tiling on $M_g$, with $g>1$, and let $i\in F$.  The centre of the algebra $A_F$ is isomorphic to a subalgebra of the algebra of $i$-loops.  It follows that, by considering lifts of paths to the universal covering space, we can consider the algebra $\Cp[\pi_1(M_g)]$ as a $Z(A_F)$ algebra, and there is a surjective $Z(A_F)$ algebra map from the algebra of $i$-loops to $\Cp[\pi_1(M_g)]$.  It follows that $A_F$ is never module finite over its centre, and there is no chance of it being a noncommutative crepant resolution of its centre.  We are also in bad shape for trying to use Proposition 7.2.14 of \cite{ginz} to establish Serre duality relative to $\mathcal{K}_{Z(A_F)}$, in fact this dualizing complex is clearly concentrated in the wrong degree.
\end{rem}

We finish with some thoughts on the Donaldson-Thomas theory of the `higher genus' 3-Calabi-Yau algebras.\smallbreak
We are fortunate, in the genus 1 case, that the fundamental group of the real torus is Abelian.  We can uniquely decompose $i$-loops (up to F-term equivalence) as $v\omega_i^t$, where $v$ is a minimal path between the endpoints of the loop, and $t$ is some nonnegative integer.  After covering $\mathbb{R}^2$ with fundamental domains of the torus, the $T_W$ weight of $u$ is given uniquely by the fundamental domain that $h(u)$ lives in.  In the non-Abelian case we are not so fortunate.  \smallbreak
Say $h(u)$ is in a fundamental domain corresponding to the element 
$\alpha\in\pi_1(M_g)$, with $g>1$.  Then any torus weight on $u$ can only be used to tell us the image of $\alpha$ under the natural noninjective map $\pi_1(M_g)\rightarrow H_1(M_g)$, the Abelianization of $\pi_1(M_g)$.  The result is that we lose our assurance that the (topological) $T_W$-fixed locus is the set of $G_1$-graded modules, and should expect to see Grassmannians, or worse, in the $T_W$-fixed locus.  Furthermore, our argument for Proposition (\ref{trivialweights}) will not do as it stands, and so more work must be done to show that the $T_W$-fixed subscheme is reduced.  We conclude that the Donaldson-Thomas picture for higher genus consistent brane tilings is more complicated, though perhaps more interesting too.

\bibliographystyle{amsplain}
\bibliography{conspreprint040410}
\end {document}

%% file: fterm.tex
\begin{picture}(0,0)%
\includegraphics{fterm.pstex}%
\end{picture}%
\setlength{\unitlength}{2901sp}%
\begingroup\makeatletter\ifx\SetFigFontNFSS\undefined%
\gdef\SetFigFontNFSS#1#2#3#4#5{%
  \reset@font\fontsize{#1}{#2pt}%
  \fontfamily{#3}\fontseries{#4}\fontshape{#5}%
  \selectfont}%
\fi\endgroup%
\begin{picture}(4139,3103)(744,-2895)
\put(1081,-1681){\makebox(0,0)[lb]{\smash{{\SetFigFontNFSS{11}{13.2}{\familydefault}{\mddefault}{\updefault}{\color[rgb]{0,0,0}u}%
}}}}
\put(4456,-1636){\makebox(0,0)[lb]{\smash{{\SetFigFontNFSS{11}{13.2}{\familydefault}{\mddefault}{\updefault}{\color[rgb]{0,0,0}v}%
}}}}
\end{picture}%

%% file: windup2.tex
\begin{picture}(0,0)%
\includegraphics{windup2.pstex}%
\end{picture}%
\setlength{\unitlength}{2901sp}%
\begingroup\makeatletter\ifx\SetFigFontNFSS\undefined%
\gdef\SetFigFontNFSS#1#2#3#4#5{%
  \reset@font\fontsize{#1}{#2pt}%
  \fontfamily{#3}\fontseries{#4}\fontshape{#5}%
  \selectfont}%
\fi\endgroup%
\begin{picture}(6404,3657)(683,-3040)
\put(6280,-1430){\makebox(0,0)[lb]{\smash{{\SetFigFontNFSS{8}{9.6}{\rmdefault}{\mddefault}{\updefault}{\color[rgb]{0,0,0}j}%
}}}}
\put(4314,-2329){\makebox(0,0)[lb]{\smash{{\SetFigFontNFSS{5}{6.0}{\familydefault}{\mddefault}{\updefault}{\color[rgb]{0,0,0}$u_1$}%
}}}}
\put(4256,-1037){\makebox(0,0)[lb]{\smash{{\SetFigFontNFSS{5}{6.0}{\familydefault}{\mddefault}{\updefault}{\color[rgb]{0,0,0}$u_2$}%
}}}}
\put(4256,-701){\makebox(0,0)[lb]{\smash{{\SetFigFontNFSS{5}{6.0}{\familydefault}{\mddefault}{\updefault}{\color[rgb]{0,0,0}$u_3$}%
}}}}
\put(3466,-1546){\makebox(0,0)[lb]{\smash{{\SetFigFontNFSS{10}{12.0}{\rmdefault}{\mddefault}{\updefault}{\color[rgb]{0,0,0}b}%
}}}}
\put(3061,-1096){\makebox(0,0)[lb]{\smash{{\SetFigFontNFSS{10}{12.0}{\rmdefault}{\mddefault}{\updefault}{\color[rgb]{0,0,0}c}%
}}}}
\put(3106,-2041){\makebox(0,0)[lb]{\smash{{\SetFigFontNFSS{10}{12.0}{\rmdefault}{\mddefault}{\updefault}{\color[rgb]{0,0,0}a}%
}}}}
\put(2341,-2806){\makebox(0,0)[lb]{\smash{{\SetFigFontNFSS{5}{6.0}{\familydefault}{\mddefault}{\updefault}{\color[rgb]{0,0,0}$u_{\infty}$}%
}}}}
\put(2431,-1546){\makebox(0,0)[lb]{\smash{{\SetFigFontNFSS{8}{9.6}{\familydefault}{\mddefault}{\updefault}{\color[rgb]{0,0,0}$i_0$}%
}}}}
\end{picture}%

%% file: newdig.tex
\begin{picture}(0,0)%
\includegraphics{newdig.pstex}%
\end{picture}%
\setlength{\unitlength}{2901sp}%
\begingroup\makeatletter\ifx\SetFigFontNFSS\undefined%
\gdef\SetFigFontNFSS#1#2#3#4#5{%
  \reset@font\fontsize{#1}{#2pt}%
  \fontfamily{#3}\fontseries{#4}\fontshape{#5}%
  \selectfont}%
\fi\endgroup%
\begin{picture}(6732,5457)(2956,-4999)
\put(5851,-2896){\makebox(0,0)[lb]{\smash{{\SetFigFontNFSS{12}{14.4}{\familydefault}{\mddefault}{\updefault}{\color[rgb]{0,0,0}$\mathcal{L}$}%
}}}}
\put(5536,-2131){\makebox(0,0)[lb]{\smash{{\SetFigFontNFSS{8}{9.6}{\familydefault}{\mddefault}{\updefault}{\color[rgb]{0,0,0}$w_0$}%
}}}}
\put(7201,-2221){\makebox(0,0)[lb]{\smash{{\SetFigFontNFSS{8}{9.6}{\familydefault}{\mddefault}{\updefault}{\color[rgb]{0,0,0}$w_1$}%
}}}}
\put(7831,-16){\makebox(0,0)[lb]{\smash{{\SetFigFontNFSS{8}{9.6}{\familydefault}{\mddefault}{\updefault}{\color[rgb]{0,0,0}$u_1$}%
}}}}
\put(5851,-1456){\makebox(0,0)[lb]{\smash{{\SetFigFontNFSS{8}{9.6}{\familydefault}{\mddefault}{\updefault}{\color[rgb]{0,0,0}$v_0$}%
}}}}
\put(7021,-1411){\makebox(0,0)[lb]{\smash{{\SetFigFontNFSS{8}{9.6}{\familydefault}{\mddefault}{\updefault}{\color[rgb]{0,0,0}$v_1$}%
}}}}
\put(5221,-4921){\makebox(0,0)[lb]{\smash{{\SetFigFontNFSS{9}{10.8}{\familydefault}{\mddefault}{\updefault}{\color[rgb]{0,0,0}$i_0$}%
}}}}
\put(4546,-286){\makebox(0,0)[lb]{\smash{{\SetFigFontNFSS{8}{9.6}{\familydefault}{\mddefault}{\updefault}{\color[rgb]{0,0,0}$u_0$}%
}}}}
\end{picture}%

%% file: expansion.tex
\begin{picture}(0,0)%
\includegraphics{expansion.pstex}%
\end{picture}%
\setlength{\unitlength}{4144sp}%
\begingroup\makeatletter\ifx\SetFigFontNFSS\undefined%
\gdef\SetFigFontNFSS#1#2#3#4#5{%
  \reset@font\fontsize{#1}{#2pt}%
  \fontfamily{#3}\fontseries{#4}\fontshape{#5}%
  \selectfont}%
\fi\endgroup%
\begin{picture}(5139,2282)(1195,-4629)
\put(3504,-3169){\makebox(0,0)[lb]{\smash{{\SetFigFontNFSS{14}{16.8}{\familydefault}{\mddefault}{\updefault}{\color[rgb]{0,0,0}$\kappa^{-1}(n+1)$}%
}}}}
\put(3903,-4294){\makebox(0,0)[lb]{\smash{{\SetFigFontNFSS{17}{20.4}{\familydefault}{\mddefault}{\updefault}{\color[rgb]{0,0,0}$\mathfrak{T}_n$}%
}}}}
\end{picture}%

%% file: consjoa020511.bbl
\providecommand{\bysame}{\leavevmode\hbox to3em{\hrulefill}\thinspace}
\providecommand{\MR}{\relax\ifhmode\unskip\space\fi MR }
\providecommand{\MRhref}[2]{%
  \href{http://www.ams.org/mathscinet-getitem?mr=#1}{#2}
}
\providecommand{\href}[2]{#2}
\begin{thebibliography}{1}

\bibitem{kazushi}
{A. Ishii, K. Ueda}, \emph{On moduli spaces of quiver representations
  associated with dimer models}, arXiv:0710.1898 (2007).

\bibitem{twistedquiverbundles}
{B. Gothen, A. King}, \emph{Homological algebra of twisted quiver bundles}, J.
  {L}ondon {M}ath. {S}oc (2005), no.~71(1), 85--99.

\bibitem{longout}
N.~Broomhead, \emph{{Dimer models and Calabi-Yau algebras}}, Ph.D. thesis,
  2009.

\bibitem{ginz}
V.~Ginzburg, \emph{{C}alabi-{Y}au algebras}, arxiv:math/0612139.

\bibitem{ObstructionsHilbert}
{K. Behrend, B. Fantechi}, \emph{Symmetric obstruction theories and {H}ilbert
  schemes of points on threefolds}, Algebra {N}umber {T}heory (2008), no.~3,
  313--345.

\bibitem{MR}
{S. Mozgovoy, M. Reineke}, \emph{On the noncommutative {D}onaldson-{T}homas
  invariants arising from brane tilings}, Arxiv:0809:0117.

\bibitem{conifold}
B.~Szendr\H{o}i, \emph{Non-commutative {D}onaldson-{T}homas theory and the
  conifold}, Geom. {T}opol. \textbf{12} (2008), no.~2, 1171--1202.

\bibitem{dualizingcomplexes}
A.~Yekutieli, \emph{Dualizing complexes, {M}orita equivalence and the derived
  {P}icard group of a ring}, J. {L}ondon {M}ath. {S}oc. \textbf{60} (1999),
  no.~3, 723--746.

\end{thebibliography}
